\documentclass[10pt,reqno]{amsart} 
\usepackage{epsfig,amssymb,amsmath,version}
\usepackage{amssymb,version,graphicx,fancybox,mathrsfs}
\usepackage{url,hyperref}
\usepackage{subfigure}
\usepackage{color}
\usepackage{stmaryrd}
\usepackage{multirow}
\usepackage{booktabs,siunitx}
\usepackage{booktabs}
\usepackage{multicol}
\usepackage[ruled,linesnumbered]{algorithm2e}
\allowdisplaybreaks
\catcode`\@=11 \theoremstyle{plain}
\@addtoreset{equation}{section}   

\@addtoreset{figure}{section}
\renewcommand\thefigure{\thesection.\@arabic\c@figure}
\newtheorem{thm}{\bf Theorem}

\newtheorem{cor}{\bf Corollary}

\newtheorem{lmm}{\bf Lemma}

\theoremstyle{remark}
\newtheorem{rem}{\bf Remark}[section]

\definecolor{ligreen}{rgb}{0.0, 0.3, 0.0}

\definecolor{darkblue}{rgb}{0.0, 0.0, 0.55}

\definecolor{anti-flashwhite}{rgb}{0.55, 0.57, 0.68}

\newcommand{\bs}[1]{\boldsymbol{#1}}

\graphicspath{{./Figures/} } 

\begin{document}
\bibliographystyle{plain}

\title[An algorithm for multiple solutions] {A spectral Levenberg-Marquardt-Deflation method for multiple solutions of semilinear elliptic systems}
\author[
	L. Li,\, Y. Zhou,\,   P. Xie\,\, $\&$\,  H. Li
	]{
		\;\; Lin Li${}^{\dag}$, \;\; Yuheng Zhou${}^{\sharp}$,\;\; 
  Pengcheng Xie${}^{\S}$ \;\; and\;\;Huiyuan Li ${}^{\ddag}$
		}

 \thanks{${}^{\dag}$School of Mathematics and Physics, University of South China, Hengyang, China. Email:   lilinmath@usc.edu.cn.\\
\indent ${}^{\sharp}$School of Mathematics and Physics, University of South China, Hengyang, China. Email: 15321363156@163.com (Y. Zhou).\\
 \indent ${}^{\S}$Corresponding author. Applied Mathematics and Computational Research Division, Lawrence Berkeley National Laboratory, 1 Cyclotron Road, Berkeley, 94720, CA, USA. Email: pxie@lbl.gov, pxie98@gmail.com (P. Xie). The work was done before the corresponding author joined LBNL. \\ 
   \indent${}^{\ddag}${State Key Laboratory of Computer Science/Laboratory of Parallel Computing, Institute of Software, Chinese Academy of Sciences, Beijing 100190, China}. Email: huiyuan@iscas.ac.cn (H. Li).}

\keywords{Multiple solutions, Legendre-Galerkin method, Levenberg-Marquardt method, deflation} \subjclass[2000]{65N35, 65N22, 65F05, 65L10}

\begin{abstract} Many nonlinear differential equations arising from practical problems may permit nontrivial multiple solutions relevant to applications, and these multiple solutions are helpful to deeply understand these practical problems and to improve some applications. Developing an efficient numerical method for finding multiple solutions is very necessary due to the nonlinearity and multiple solutions of these equations. Moreover, providing an efficient iteration plays an important role in successfully obtaining multiple solutions with fast and stable convergence. In the current paper, an efficient algorithm for finding multiple solutions of semilinear elliptic systems is proposed, where the trust region Levenberg-Marquardt method is firstly used to iterate the resulted nonlinear algebraic system. When the nonlinear term in these equations has only the first derivative, our algorithm can efficiently find multiple solutions as well. Several numerical experiments are tested to show the efficiency of our algorithm, and some solutions which have not been shown in the literature are also found and shown.
\end{abstract}
\maketitle

\section{Introduction}\label{sect1}

Many nonlinear differential equations arising from practical problems may permit nontrivial multiple solutions, which often occurs in physics, mechanics, biology, energy and engineering and so on. The in-depth study of multiple solutions is helpful to improve the understanding and application of these problems \cite{2002A, davis1960introduction, Tadmor2012A}. Unfortunately, most equations do not have explicit solutions. Moreover, the PDE theory only provides some powerful analytic solution techniques for special cases (such as radial symmetrical case). The development and study of efficient numerical methods for finding multiple solutions become very meaningful, and are attracting the attention of many researches around the world.

Recently, in \cite{2022Two}, we proposed a spectral trust-region deflation method for finding multiple solutions of a single nonlinear equation. While for some coupled nonlinear equations (or models) with multiple solutions, the spectral trust-region deflation method presented in \cite{2022Two} doesn't seem to be very effective without any further improvement. To be specific, we identify multiple solutions of the following problem
\begin{equation}\label{New1.1}
-\Delta \vec{\textit{\textbf{u}}} = \vec{G}(\textit{\textbf{x}}, \vec{\textit{\textbf{u}}})   \quad\quad     \textit{\textbf{x}}\in \Omega
\end{equation}
supplemented with some boundary conditions, where $\Omega$ is a bounded domain in $R^{d}$, $d = 1, 2, \cdots$, and $\vec{G}$ is a nonlinear vector function of $\vec{\textit{\textbf{u}}}$. Comparing with finding multiple solutions of a single nonlinear equation, finding multiple solutions of \eqref{New1.1} efficiently becomes more difficult and challenging due to the significant increase in the number of discrete unknown variables. For example, when we consider the case $d = 2, \vec{\textit{\textbf{u}}} := (u, v)$ with homogeneous Dirichlet boundary conditions, and assume that the numerical solution (denoted by $u_{N}(x, y)$ and $v_{N}(x, y)$) is expanded as
\begin{equation}\label{NNew1.1}
u_{N}(x, y) = \sum_{i,j=0}^{N}\tilde{u}_{ij}\phi_{i}(x)\phi_{j}(y) \quad \textrm{and} \quad  v_{N}(x, y) = \sum_{i,j=0}^{N}\tilde{v}_{ij}\phi_{i}(x)\phi_{j}(y),
\end{equation}
where $\phi_{i}(x) = L_{i}(x) - L_{i+2}(x)$, $L_{i}(x) (0\leq i \leq N + 2)$ are Legendre polynomials, and $\{\tilde{u}_{ij}\}$ and $\{\tilde{v}_{ij}\}$ are expanding coefficients to be solved. From \eqref{NNew1.1}, we need to solve $2(N+1)^2$ unknown variables. The computational complexity increases dramatically with increasing $N$, which leads to that finding multiple solutions of \eqref{New1.1} (or solving these unknown variables in \eqref{NNew1.1}) becomes more and more difficult. On the other hand, in some coupled PDEs with multiple solutions, such as two types of noncooperative systems in \cite{2010A}, $\vec{G} \in C^{1}[-1, 1]$ in \eqref{New1.1} is only satisfied with fixed parameters. In other words, the algorithms using second derivative information can't be used to find multiple solutions of \eqref{New1.1} in the situation, which leads to that the spectral trust-region deflation method given in \cite{2022Two} can't be applied directly to \eqref{New1.1} due to necessary second derivative information (i.e. the Hessian matrix). Based on these difficulties above, finding multiple solutions of \eqref{New1.1} efficiently is more challenging. Therefore, in this paper, our aim is to improve (or extend) the spectral trust-region deflation method presented in \cite{2022Two}, and to design a new numerical method for finding multiple solutions of \eqref{New1.1} efficiently.

Before we introduce our approach, it is necessary to elaborate on the relevant background and motivations. In 1993, Choi and Mckenna \cite{choi1993mountain} proposed the mountain pass algorithm (MPA) for multiple solutions of semilinear elliptic problems. Subsequently, in \cite{xie2005improved}, Xie at al. pointed out that the MPA was feasible for finding two solutions of mountain pass type with Morse index 1 or 0. However, As said in \cite{ding1999high}, the MPA may fail to locate the sign-changing solutions, where Ding et al. proposed a high linking algorithm (HLA) for such solutions.
In 2001, Li and Zhou \cite{li2001minimax} proposed the minimax algorithm (MNA) for finding multiple solutions of nonlinear equations, and more recent advancements are presented in \cite{2005A, 2007numerical, 2008numerical, 2005Instability}. It is worth pointing out that all of these methods above require that nonlinear differential equations have the variational structure, where the variational structure plays an important role in designing the algorithm to find multiple solutions. However, many differential equations with multiple solutions have no the variational structure, which makes the methods above are not applicable. This also leads to the second category of the existing methods for multiple solutions, i.e. some numerical methods (e.g. spectral method or finite difference) are used to discretize differential equations with multiple solutions. Then some iterative methods are provided to find multiple solutions of the resulting nonlinear algebraic system (NLAS). Along this line, the search-extension method (SEM) was proposed by Xie at al \cite{chen2004search}, and some improvements to SEM have been found in \cite{xie2006improved, xie2015augmented}. In \cite{allgower2006solution, allgower2009application},
with the finite difference discretization, Allgower et al. proposed the homotopy continuation method for finding multiple solutions of the NLAS, and some interesting recent works have been inspired by this method, e.g. \cite{hao2014bootstrapping, 2018Two, zhang2013eigenfunction}. It is worth pointing out that in these methods, the Newtonian method is often chosen to iterate the NLAS with an initial guess.  However, as we known, the main disadvantage of the Newtonian method is that the iteration is sensitive to the initial guess and the inverse of the Jacobian matrix should exist. While in \cite{2022Two}, we firstly used the trust-region method to replace the Newtonian method for finding multiple solutions, and the numerical results concerning computational efficiency have also been greatly improved. But when the trust-region method is continued to iterate the NLAS from \eqref{New1.1}, the corresponding computational efficiency will be reduced due to the expensive computational cost of the Hessian matrix. More importantly, if the Hessian matrix doesn't exist, the trust-region method given in \cite{2022Two} can't be directly used for  iterating the NLAS. So in current work, combing with the deflation technique, we will design a more efficient iteration and substantially improve the efficiency for finding multiple solutions of \eqref{New1.1}. To guarantee a good approximation of the resulted NLAS to the original differential equations with a relatively low computational cost, the Legendre Spectral-Galerkin method is mainly used to discretize \eqref{New1.1}. With only first derivative information, the trust region Levenberg-Marquardt method is firstly introduced to iterate the resulted NLAS and the deflated system. Here our method is dubbed as the spectral trust region LM-Deflation method. Compared with the existing methods, the main differences and advantages of the spectral trust region LM-Deflation method reside two respects:
\begin{itemize}
  \item When the number of unknown variables required to be solved increases significantly, the spectral trust region LM-Deflation method can reduce efficiently the computational complexity in the whole iterative process through ignoring the Hessian matrix. Meanwhile, the method can also be used to solve \eqref{New1.1} with $\vec{G} \in C^{1}[a, b]$.
  \item Comparing with several existing methods, the method is capable of obtaining multiple solutions with fast and stable convergence, which allows us to start even with the same initial guesses for multiple solutions based on the deflation technique.
\end{itemize}

The remainder of this paper is organized as follows. In section \ref{sect2}, we describe the spectral
trust region LM-Deflation method. In section \ref{sect3}, numerical experiments are provided to demonstrate the efficiency of the method. In section \ref{sect4}, we conclude the paper with some remarks.

\section{An algorithm for computing multiple solutions}\label{sect2}

\begin{figure}[!ht]
\begin{centering}
\includegraphics[width=6.5cm,height=5.5cm]{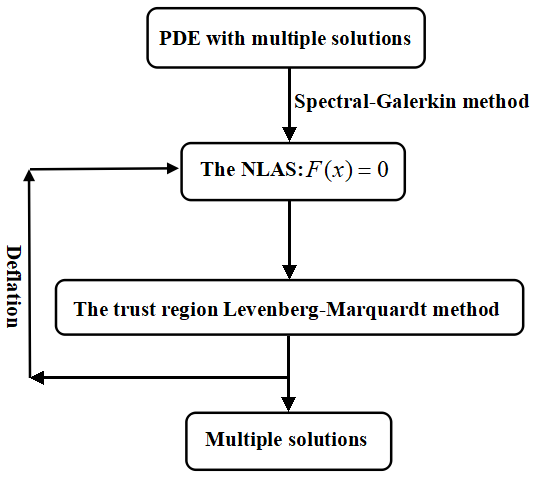}\vspace{-0.25cm}
\caption{The simple flow configuration of our algorithm.}\label{NFg2}
\end{centering}
\end{figure}

The main purpose of this section is to propose an efficient algorithm for computing multiple solutions of \eqref{New1.1}. The simple flow configuration of the algorithm is presented in Fig.\ref{NFg2}, where the Spectral-Galerkin method is used to discretize \eqref{New1.1}, the trust region Levenberg-Marquardt method is designed to iterate the resulted NLAS and the deflated system, and the deflation technique is used to deflate a known solution and then obtain a another solution. For added clarity, the forthcoming section is divided into four parts: the first part - the Spectral-Galerkin method to \eqref{New1.1}, the second part - the trust region Levenberg-Marquardt method for iterating the resulted NLAS and the deflated system, the third part - the deflation technique for finding multiple solutions and the fourth part - a summary of our algorithm.

\subsection{The Spectral-Galerkin method}\label{sec2.1}
We consider the following general system of partial differential equations
\begin{equation}\label{modify08121}
\mathcal{L}{\bs u} + \mathcal{N}{\bs u} = \mathcal{{\bs 0}},  \quad\quad  {\bs x}\in \Omega,
\end{equation}
where ${\bs u}:= (u_{1}({\bs x}), u_{2}({\bs x}), \cdots, u_{n}({\bs x}))^{T}$ is a vector function of ${\bs x}$, $\mathcal{L}$ and $\mathcal{N}$ are linear and nonlinear operator, respectively. Here the linear operator $\mathcal{L}$ may be $-\Delta$ or other operators, and in this paper we mainly focus on the former. The nonlinear operator $\mathcal{N} := (\mathcal{N}_{1}({\bs u}), \mathcal{N}_{2}({\bs u}), \cdots, \mathcal{N}_{n}({\bs u}))^{T}$ is defined. When the basis functions $\{\phi_{k}\}_{k=0}^{N}$ satisfied corresponding boundary conditions are chosen, the resulting $\{u_{i}({\bs x})\}_{i=1}^{n}$ can be expanded as follows
\begin{equation}\label{modify08122}
u_{i}({\bs x}) = \sum_{j=0}^{N}a^{i}_{j}\phi_{j}({\bs x}),     \quad\quad   1 \leq i \leq n.
\end{equation}
Based on the Spectral-Galerkin method, substituting \eqref{modify08122} into \eqref{modify08121} yields
the resulting nonlinear algebraic system, i.e.,
\begin{equation}\label{modify08123}
\mathcal{A}{\bs a} + \mathcal{F}({\bs a}) = {\bs 0},
\end{equation}
where
\begin{equation*}
{\bs a} = (a^{1}_{0},\, a^{1}_{1},\, \cdots, \,a^{1}_{N},\, a^{2}_{0},\, \cdots,\, a^{2}_{N},\, \cdots, \, a^{n}_{0},\, \cdots,\, a^{n}_{N})^{T},
\end{equation*}
\begin{equation*}
[{\bs S}]_{ij} = -\int_{\Omega}\Delta \phi_{j}({\bs x})\phi_{i}({\bs x})d{\bs x},  \quad\quad
\mathcal{A} = \begin{bmatrix}
{\bs S}\\
&\ddots &  &\text{{\huge 0}}\\
&  & \ddots\\
& \text{{\huge 0}}  &  & \ddots\\
& & & & {\bs S}
\end{bmatrix},
\end{equation*}
\begin{equation*}
\mathcal{F}({\bs a}) = (\int_{\Omega}\mathcal{N}_{1}({\bs u})\phi_{0}({\bs x})d{\bs x},\cdots, \int_{\Omega}\mathcal{N}_{1}({\bs u})\phi_{N}({\bs x})d{\bs x}, \int_{\Omega}\mathcal{N}_{2}({\bs u})\phi_{0}({\bs x})d{\bs x}, \cdots, \int_{\Omega}\mathcal{N}_{n}({\bs u})\phi_{N}({\bs x})d{\bs x})^{T}.
\end{equation*}
Here it is worth pointing out that ${\bs a}$ in \eqref{modify08123} is the unknown vector to be solved.
To increase the clarity of the above process, as an illustrative example, we consider \eqref{modify08121} with homogeneous Dirichlet boundary conditions, $d = 2$ (i.e. $\Omega = (x, y)$), ${\bs u} = (u_{1}(x, y), u_{2}(x, y))$ and $\mathcal{N} = (\mathcal{N}_{1}({\bs u}), \mathcal{N}_{2}({\bs u}))$, and the Legendre-Galerkin method is used to discretize it. Let ${\mathbb P}_N$ be the set of the polynomials of degree at most $N,$ and ${\mathbb P}_N^0=\{\phi\in {\mathbb P}_N\,:\, \phi(\pm 1)=0\}.$ The Legendre-Galerkin approximation is to find $(u^{1}_N, u^{2}_{N})\in X_N^0=({\mathbb P}_N^0)^2$ such that
 \begin{equation}\label{New2.3B}
 \begin{cases}
(\nabla u^{1}_N, \nabla \phi_N) = (I_N \mathcal{N}_{1}, \phi_{N}),\\
(\nabla u^{2}_N, \nabla \phi_N) = (I_N \mathcal{N}_{2}, \phi_{N}), \\
\end{cases}
\end{equation}
for $\forall \phi_{N}\in {\mathbb P}_N^0$, where $I_N$ is the Legendre-Gauss-Lobatto tensorial interpolation operator with $N+1$ points in each coordinate direction. Based on the property of the Legendre polynomial \cite{shen2011spectral}, we introduce the basis of ${\mathbb P}^0_{N}$ as 
\begin{equation*}
  \phi_{k}(x) = L_{k}(x) - L_{k+2}(x),\quad 0\le k\le N,
  \end{equation*}
and write
\begin{equation}\label{2.5}
\begin{split}
& u^{1}_{N} = \sum^{N}_{k,j = 0}a^{1}_{kj}\phi_{k}(x)\phi_{j}(y), \quad  u^{2}_{N} = \sum^{N}_{k,j = 0}a^{2}_{kj}\phi_{k}(x)\phi_{j}(y),   \\
& \hat{a}_{kj} = \int_{I}\phi'_{j}(x)\phi'_{k}(x)dx,   \quad\quad  A = (\hat{a}_{kj})_{k, j = 0, 1, \cdots, N};\notag\\
& b_{kj} = \int_{I}\phi_{j}(x)\phi_{k}(x)dx,   \quad\quad  B = (b_{kj})_{k, j = 0, 1, \cdots, N};\notag\\
& \mathcal{F}_{1, kj} = (I_{N}\mathcal{N}_{1}, \phi_{k}(x)\phi_{j}(y))    \quad  \;  \mathcal{F}_{2, kj} = (I_{N}\mathcal{N}_{2}, \phi_{k}(x)\phi_{j}(y)).\notag\\
\end{split}
\end{equation}
As a result, ${\bs S}, {\bs a}$ and $\mathcal{F}$ in \eqref{modify08123} become
\begin{equation*}
{\bs S} = A \otimes B + B \otimes A^{T},
\end{equation*}
\begin{equation*}
{\bs a} = (a^{1}_{00}, a^{1}_{10}, \cdots, a^{1}_{N0}, a^{1}_{01}, \cdots, a^{1}_{N1}, \cdots, a^{1}_{NN}, a^{2}_{00}, \cdots, a^{2}_{NN})^{T},
\end{equation*}
\begin{equation*}
\mathcal{F} = (\mathcal{F}_{1, 00}, \mathcal{F}_{1, 10}, \cdots, \mathcal{F}_{1, N0}, \cdots, \mathcal{F}_{1, NN},  \mathcal{F}_{2, 00}, \mathcal{F}_{2, 10}, \cdots, \mathcal{F}_{2, N0}, \cdots, \mathcal{F}_{2, NN})^{T}.
\end{equation*}
Here $\otimes$ denotes the tensor product operator, i.e. $A \otimes B = (Ab_{ij})_{i,j = 0, 1, \cdots, N}$. Now \eqref{modify08123} becomes a nonlinear system due to $\mathcal{F}_{1, kj}$ and $\mathcal{F}_{2, kj}$, which will be solved by the trust region Levenberg-Marquardt method. Here it is worth pointing out that in the computational process, we only need to evaluate $\mathcal{F}_{1, kj}$ and $\mathcal{F}_{2, kj}$ for given ${\bs u}$ that can be implemented efficiently by the pseduo-spectral technique described in \cite[Ch.\! 4]{shen2011spectral} in terms of coefficients.

\subsection{The trust region Levenberg-Marquardt method}\label{sect2.3}
Formally, based on the legendre-Galerkin method to \eqref{modify08121}, we can obtain a nonlinear algebraic system like \eqref{modify08123}, and write this system in a more elegant form, i.e.,
\begin{equation}\label{4.1}
{\bs F} =\big({F}_{1}, \; {F}_{2}, \; \cdots, \; {F}_{n}\big)^{T}
\end{equation}
with the unknown vector ${\bs a} = (a_{1},\, \cdots,\, a_{n})^{T}$, where ${F}_{i}\, (i = 1, \cdots, n)$ are functions of the unknown vector ${\bs a}$. Our aim is to find ${\bs a}$ to satisfy ${\bs F}({\bs a}) \equiv {\bs 0}$ using an iterative method. Along the line in \cite{2022Two}, we firstly reformulate the zero-finding problem \eqref{4.1} as the optimisation problem using the nonlinear least-squares method:
\begin{equation}
\min_{\bs a \in {\mathbb R}^{n}}Q(\textbf{\textit{a}}),\quad Q(\bs a):= \frac{1}{2}\big\|\bs{F}(\bs a)\big\|^{2}_{2} = \frac{1}{2}\sum_{i = 1}^{n}{F}^2_{i}(\textbf{\textit{a}}). \label{4.2}
\end{equation}
Before we present a method to solve \eqref{4.2}, the Jacobian matrix, the gradient and Hessian matrix should be introduced as follows:
\begin{equation*}\label{New4.4}
\bs J(\bs a) = {\bs F}'(\textbf{\textit{a}}) = (\nabla {F}_{1}(\textbf{\textit{a}}), \nabla {F}_{2}(\textbf{\textit{a}}), \cdots, \nabla F_{n}(\textbf{\textit{a}}))^{T},
\end{equation*}
\begin{equation*}
{\bs g}({\bs a}) = \nabla Q({\bs a}) = \bs{J}^{T}(\bs{a}) {\bs F(\bs {a})}   \quad \textrm{and} \quad   {\bs G}({\bs a}) = \nabla^2 Q({\bs a}) = {\bs J}^{T}(\bs{a}) \bs J({\bs a}) + \bs S({\bs a}),
\end{equation*}
where ${\bs S(\textbf{\textit{a}})} = \sum_{i = 1}^{n}F_{i}(\textbf{\textit{a}})\nabla^2 F_{i}(\textbf{\textit{a}})$. so far, in term of updating the unknown vector ${\bs a}$ in \eqref{4.2}, there may exist three formats as follows:
\begin{equation}\label{modify08141}
\begin{split}
& (\textbf{\textrm{I}})\quad\quad\;\;\; {\bs a}_{k+1} = {\bs a}_{k} - {\bs G}({\bs a}_{k})^{-1}{\bs g}({\bs a}_{k}),  \quad\quad  (\textrm{Gauss-Newtonian method})\\
& (\textbf{\textrm{II}})\quad\quad\;\, {\bs a}_{k+1} = {\bs a}_{k} - (J({\bs a}_{k})^{T}J({\bs a}_{k}))^{-1}{\bs g}({\bs a}_{k}),\\
&  (\textbf{\textrm{III}})\quad\quad {\bs a}_{k+1} = {\bs a}_{k} - (J({\bs a}_{k})^{T}J({\bs a}_{k}) + \mu_{k}I)^{-1}{\bs g}({\bs a}_{k}),  \quad   (\textrm{Levenberg-Marquardt method})\\
\end{split}
\end{equation}

\noindent where $\mu_{k} \geq 0$ is a parameter being updated from iteration to iteration.
Obviously, the main disadvantage of the Gauss-Newtonian method is that the second-order term ${\bs S({\bs a})}$ in the Hessian matrix ${\bs G({\bs a})}$ is difficult or expensive to compute. Furtherly, as said in the introduction above, the Hessian matrix ${\bs G({\bs a})}$ doesn't exist in some situations, which means that the Gauss-Newtonian method can't be used to solve \eqref{4.2}. On the other hand, when ignoring ${\bs S({\bs a})}$, we use the format $\textrm{\textbf{II}}$ to solve \eqref{4.2}, where it is also required that Jacobian matrix ${\bs J}({\bs a})$ has full column rank. In other words, if ${\bs J}({\bs a})$ is rank-deficient, then the format $\textrm{\textbf{II}}$ cannot work well. In addition, based on our experience from \cite{2022Two}, we have realized that computing a full Hessian matrix can reduce computational efficiency. To overcome these difficulties, here we firstly use the Levenberg-Marquardt method to solve \eqref{4.2} in the process of seeking multiple solutions. Let
\begin{equation}\label{4.7}
{\bs s}_{k} := {\bs a}_{k+1} - {\bs a}_{k} = -({\bs J}({\bs a}_k)^{T}{\bs J}({\bs a}_k) + \mu_k I)^{-1}{\bs g}({\bs a}_k),
\end{equation}
then it is easy to verify that the Levenberg-Marquardt step ${\bs s}_{k}$ in \eqref{4.7} is the solution of the following optimization problem
\begin{equation}
\min_{{\bs s}\in R^{n}}\|{\bs F}({\bs a}_{k}) + {\bs J}({\bs a}_{k}){\bs s}\|^2_2 + \mu_{k}\|{\bs s}\|_2^2.
\end{equation}
Meanwhile, the Levenberg-Marquardt method is also one of the trust-region methods. To be specific, the following trust-region model is considered:
\begin{equation}\label{4.6}
\begin{split}
& \min_{\bs s\in {\mathbb B}_{h_k}} q^{(k)}({\bs s}) := Q({\bs a}_{k}) + {\bs g}({\bs a}_{k})^{T}{\bs s} + \frac{1}{2}{\bs s}^{T}{\bs J}({\bs a}_{k})^{T}{\bs J}({\bs a}_{k}){\bs s},\quad
\end{split}
\end{equation}
where $\Delta_{k}$ represents the trust-region radius and the trust region  ${\mathbb B}_{h_k}:=\{\bs s := \bs x - \bs x_{k}\in {\mathbb R}^n\,:\,\|\bs s\|\le \Delta_{k} = {\bs s}_k\}$\cite{2012ZZhang,2014Zhangzaikun,xie2023least,xie2023derivative,xie2024derivative,xie2024newtwodimensionalmodelbasedsubspace,xie2025remuregionalminimalupdating}. Then we can derive that ${\bs s}_k$ in \eqref{4.7} is the global optimal solution of \eqref{4.6} \cite{sun2006optimization,XIE2025116146,}. Here it is worth pointing out that the main difference between the Levenberg-Marquardt method (i.e. \eqref{4.6}) and other trust region methods is that other trust region methods update the trust region radius $\Delta_{k}$ directly, while the Levenberg-Marquardt method updates the parameter $\mu_{k}$, which in turn modifies the value $\Delta_{k}$ implicitly. Moreover, we can see that $(J({\bs a}_{k})^{T}J({\bs a}_{k}) + \mu_{k}I)$ in \eqref{modify08141} is also positive define, indicating that its inverse can be well defined. When $\mu_{k}$ in \eqref{modify08141} is changed, the Levenberg-Marquardt method allows choosing any direction between Gauss-Newtonian direction and the steepest descent direction, i.e., when $\mu_{k} = 0$, it reduces to the Gauss-Newtonian direction. While $\mu_{k}$ is large enough, the produced direction is close to the steepest descent direction. Furtherly, in \cite{sun2006optimization}, Yuan et al. have described the relation between ${\bs s}$ and $\mu$, i.e., let ${\bs s} = {\bs s}(\mu)$ is a solution of
\begin{equation}
(J({\bs a})^{T}J({\bs a}) + \mu I) {\bs s} = - {\bs g}({\bs a}),
\end{equation}
then we can conclude the following theorem, i.e.,
\begin{thm}\label{thm26}  When $\mu$ increases monotonically from zero, $\|{\bs s}(\mu)\|_2$ in \eqref{4.6} will decrease strictly monotonically.
\end{thm}
To obtain a proper trust-region radius $\Delta_{k}$ for determining ${\bs s}_{k}$ and $\mu_{k}$, we define
\begin{equation}\label{rkval}
r_{k} = \frac{Q({\bs a}_{k}) - Q({\bs a}_{k+1})}{q^{(k)}(\bs 0) - q^{(k)}({\bs s}_{k})}
\end{equation}
for $q^{(k)}({\bs s}_{k})$ and the objective function value $Q({\bs a}_{k+1})$. \eqref{rkval} is viewed as an indicator for the expansion and contraction the trust region $\Delta_{k}$. If $r_{k}$ is close to 1, it means there is good agreement, we can expand the trust-region radius for the next iteration. Based on the relation between ${\bs s}_{k}$ and $\mu_{k}$ in theorem \ref{thm26}, we should reduce $\mu_{k+1}$ for this case. In current work we choose $\mu_{k+1} = 0.1\mu_{k}$. If $r_{k}$ is close to zero or negative, the trust-region radius should be shrinked. The corresponding value of $\mu_{k+1}$ should be increased, and $\mu_{k+1} = 10\mu_{k}$ is chosen. Otherwise, we do not alter the trust-region radius.

Next, we will discuss research results on the convergence of the Levenberg-Marquardt method. With some conditions, the Levenberg-Marquardt method enjoys the desirable convergence with a local superlinear rate of convergence. For the convenience of subsequent discussion, we firstly introduce the definition of the local error bound.

\noindent\textbf{Definition of the local error bound.} Let $X^{*}$ be the solution set of \eqref{4.1} and ${\bs a}^{*} \in X^{*}$, if there exists a positive constant $c > 0$ such that
\begin{equation}
\|{\bs F}({\bs a})\|_2 > c\, \textrm{dist}({\bs a}, X^{*}), \quad\quad \forall {\bs a}\in N({\bs a}^{*}, b) = \{{\bs a}|\, \|{\bs a} - {\bs a}^{*}\|_2 \leq b\} \cap X^{*} \neq \emptyset,
\end{equation}
where
\begin{equation*}
\textrm{dist}({\bs a}, X^{*}) = \inf_{{\bs y}\in X^{*}}\|{\bs a} - {\bs y}\|_{2}.
\end{equation*}
Then, we call that ${\bs F}({\bs a})$ is a local error bound on some neighbourhood of ${\bs a}^{*} \in X^{*}$.
As said in \cite{2012Fan}, if the Jacobian matrix is nonsingular at the solution ${\bs a}^{*}$ of \eqref{4.1} and if the initial guess is chosen sufficiently close to ${\bs a}^{*}$, then the Levenberg-Marquardt method has a quadratic convergence. To improve the requirement of the nonsingularity, Yamashita and Fukushima \cite{2001Yamashita} show that under the local error bound condition, if $\mu_{k} = \|{\bs F}({\bs a}_{k})\|^2_2$ and if the initial guess is chosen sufficiently close to the solution set $X^{*}$, the Levenberg-Marquardt method converges quadratically to the solution set $X^{*}$. In \cite{2005Fan}, Fan et al. considered $\mu_{k} = \|{\bs F}({\bs a}_{k})\|^{\delta}_2$ with $\delta \in [1, 2]$ and proved that the Levenberg-Marquardt method still achieves the quadratic convergence under the same conditions. In \cite{2012Fan}, the cubic convergence of the modified Levenberg-Marquardt method was also proved under the local error bound condition. Finally, some progressive and interesting studies can be seen in \cite{2019Fan}.

\subsection{The deflation technique for finding multiple solutions.}\label{sect2.4}

In fact, the deflation technique is used to find multiple solutions of \eqref{4.1}. Its main idea is to successively modify \eqref{4.1} under consideration to eliminate known solutions, discovering additional distinct solutions. we now state the deflation process for computing multiple solutions to the nonlinear algebraic system \eqref{4.1}.

\begin{enumerate}
\setlength{\itemindent}{0.5em}
  \item[\textbf{\quad Step 1.}] Give a initial guess ${\bs a}^{0} \in R^{n}$;
  \item[\textbf{Step 2.}] Based on the Levenberg-Marquardt method, a solution of \eqref{4.1} is computed with the initial
      \item[] guess ${\bs a}^{0}$, and the solution is denoted by ${\bs r_1}$;
  \item[\textbf{Step 3.}] A deflation operator is defined as 
   \begin{equation}\label{modifydeflation1}
\Psi({\bs a}; {\bs r}_1) = (\frac{1}{\|{\bs a}-{\bs r}_1\|^{p}_{U}} + \alpha){\mathbb I},
\end{equation}
    \item[] where ${\mathbb I}$ is the identity operator, $\alpha \geq 0$ is a shift scalar, $p \in [1, \infty)$ is the deflation exponent.
    \item[] The norm $\|\cdot\|_{U}$ should be chosen for the solution space. A deflated system $\hat{{\bs F}}({\bs a})$ can be
    \item[] formed by applying \eqref{modifydeflation1} to \eqref{4.1} as
\begin{equation}
\hat{{\bs F}}({\bs a}) = \Psi({\bs a}; {\bs r}_1){\bs F}({\bs a});
\end{equation}
\item[\textbf{Step 4.}] With the same initial guess ${\bs a}^{0}$, the deflated system $\hat{{\bs F}}({\bs a})$ is solved using the Levenberg-
    \item[] Marquardt method, and the resulting solution is denoted by ${\bs r}_2$;
\item[\textbf{Step 5.}] A multiple deflation operator needs to be defined as follow:
\begin{equation}\label{NNNew2.3.1}
\Psi({\bs a}; {\bs r_1}, {\bs r_2}) = \prod \limits_{i=1}^{2}\Psi({\bs a}; {\bs r}_i),
\end{equation}
\item[] and a new deflated system $\hat{\hat{{\bs F}}}({\bs a})$ can be obtained: 
\begin{equation}
\hat{\hat{{\bs F}}}({\bs a}) = \Psi({\bs a}; {\bs r_1}, {\bs r_2}){\bs F}({\bs a});
\end{equation}
\item[\textbf{Step 6.}] Similar to \textbf{Step 4}, and the resulting solution is denoted by ${\bs r}_3$;
\item[\textbf{Step 7.}] Continue the cycle from Step 5 to Step 6, multiple solutions can be obtained, and these
    \item[] solutions are denoted by $X^{*} = \{{\bs r}_1, {\bs r}_2, \cdots\}$.
\end{enumerate}
\begin{rem} To have a better understanding and application of the deflation, let us make some
remarks.
\end{rem}

In Step 3, the main purpose of introducing the deflation operator \eqref{modifydeflation1} is to find another distinct solution. To be specific, based on the results in \cite{1971Deflation, farrell2015deflation, 2022Two}, the deflated system $\hat{{\bs F}}({\bs a})$ satisfies two properties: (i) The preservation of solutions of $\hat{{\bs F}}({\bs a})$ should be hold, i.e. for ${\bs a} \neq {\bs r}_1$, $\hat{{\bs F}}({\bs a}) = {\bf 0}$ iff ${\bs F}({\bs a}) = {\bf 0}$; (ii) The Levenberg-Marquardt method applied to $\hat{{\bs F}}({\bs a})$ will not find ${\bs r}_1$ again due to the deflation operator, i.e.,
\begin{equation*}
\lim_{{\bs a} \rightarrow {\bs r}_1}\inf \|\hat{{\bs F}}({\bs a})\|_{U} > 0,
\end{equation*}
indicating that with the deflation operator, another distinct solution can be obtained, even with the same initial guess. In addition, for simplicity, in current paper $p = 2$ and $\alpha = 1$ in \eqref{modifydeflation1} are chosen.

In Step 4, it is worth pointing out that when the same initial guess is chosen, the Levenberg-Marquardt method combing with the deflation may diverge in some cases. Motivated by \cite{2022Two}, we consider the choice of a random deviation from the obtained solution as the initial guess turns out sufficient to find all other solutions.

In Step 5, it is worth mentioning that when we consider more than two solutions, the multiple deflation operator is constructed, which may increase the computational complexity in the nonlinear iteration. The single deflation operator with the initial guess presented in Step 4 may be considered to replace the multiple deflation operator, which is helpful to reduce the computational complexity.

\subsection{Summary of the spectral Levenberg-Marquardt-Deflation algorithm}\label{sect2.5}
In summary, based on the methods above, a new and novel algorithm of computing multiple solutions of \eqref{New1.1} can be designed, which is presented in the following {\bf Algorithm 1}. \vspace{0.25cm}

\begin{center}
\begin{tabular}{l}
\toprule[0.8pt]
\textbf{Algorithm 1} - An algorithm of computing multiple solutions of \eqref{New1.1}\\\toprule[0.8pt]
\quad \textbf{Input:} $0 < \epsilon \ll 1$, $0 < \delta_{1} < \delta_{2} < 1$, $\mu_{0} = 0.01$, and initial solution set $S$ $\Leftarrow$ the empty set $\Phi$. \\
\quad \textbf{Output:} $S$ \\
\, 1: \quad Given nonlinear discrete system ${\bs F}({\bs a})$ by the Legendre-Galerkin method in the section \ref{sec2.1}.\\
\, 2: \quad{\bf {while}} multiple solutions not reached {\bf {do}}\\
\, 3: \quad\quad {\bf {while}} first-order optimality threshold or failure criterion not met {\bf {do}}\\ 
\, 4: \quad\quad\quad Compute ${\bs J}({\bs a}_{k})$ and ${\bs g}({\bs a}_{k})$; \quad $\triangleright$ Levenberg-Marquardt iteration in section \ref{sect2.3}.\\
\, 5: \quad\quad\quad If $\|{\bs g}({\bs a}_{k})\|<\epsilon$ and $|Q({\bs a}_{k})|^{1/2}<\epsilon$, stop;\\
\, 6: \quad\quad\quad Approximately solve ${\bs s}_{k}$ by \eqref{4.7};\\
\, 7: \quad\quad\quad Compute $r_{k}$ by \eqref{rkval};\\
\, 8: \quad\quad\quad If $r_{k}\geq \delta_{1}$, then ${\bs a}_{k+1} = {\bs a}_{k} + {\bs s}_{k}$; Otherwise, set ${\bs a}_{k+1} = {\bs a}_{k}$;\\
\, 9: \quad\quad\quad If $r_{k}<\delta_{1}$, then $\mu_{k+1}:=10\mu_{k}$;\\
\quad\quad\quad\quad\; If $\delta_{1}\leq r_{k} \leq \delta_{2}$, then $\mu_{k+1}:=\mu_{k}$;\\
\quad\quad\quad\quad\;  Otherwise, set $\mu_{k+1}:=0.1 \mu_{k}$;\\
 10: \quad\quad\quad Update $q^{(k)}$, set $k:= k+1$, go to step 4; \\
 11: \quad\quad {\color{blue}{end}}\\
 12: \quad\quad\quad \textbf{if} Convergence threshold met \textbf{then}\\
 13: \quad\quad\quad Add converged solution ${\bs r}$ to $S$;\\
 14: \quad\quad\quad $\hat{{\bs F}}({\bs a})$ $\leftarrow \Psi({\bs a}; {\bs r}){\bs F}({\bs a})$ \quad  $\triangleright$ the deflation technique presented in section \ref{sect2.4}.\\
 15: \quad\quad\quad Go back to step 3; \\
 16: \quad\quad\quad $\textbf{end}$\\
 17:\quad\; {\bf {end}}\\
 18: return $S$\\
\toprule[0.8pt]
\end{tabular}\vspace{0.1cm}
\end{center}

Some additional remarks about the implementation of the {\bf Algorithm 1} are listed, and the corresponding numerical tests will be proceeded in the forthcoming section.
\begin{itemize}
 \item As highlighted previously, the Levenberg-Marquardt method is firstly introduced.
   Comparing with the Newtonian method, the Levenberg-Marquardt method enables us to overcome the sensitivity of initial guesses, and substantially improve the computational efficiency.
 \item When increasing $N$, the matrix-vector multiplication can also be used to save the storage space in a computer.
\end{itemize}

\section{Numerical experiments}\label{sect3}

In this section, with several numerical experiments, we focus on the efficiency of {\bf Algorithm 1} for finding multiple solutions. One and two dimensional numerical experiments are shown in sections \ref{sect3.1}-\ref{sect3.2}, respectively. We choose $\delta_1 = 0.25, \delta_2 = 0.75,$ and $\epsilon = 10^{-13}$ in {\bf Algorithm 1}, and all programs are carried out on a server with Intel(R) Core(TM) i7-7500U (2.90 GHz) and 120GB RAM.

\subsection{ODE Examples}\label{sect3.1} We consider the following problem:
\begin{equation}\label{Nexample3.1}
\begin{cases}
d_1\frac{d^2{u}}{d{x}^2} = g_1(u, v)  \vspace{0.15cm}\\
d_2\frac{d^2{v}}{d{x}^2} = g_2(u, v)
\end{cases}  \quad \quad  \textrm{on}\; \Omega = (0, 1)
\end{equation}
with the no-flux boundary conditions: $\frac{d{u}}{d{x}} = 0$ and $\frac{d{v}}{d{x}} = 0$ in $\partial\Omega$, where $d_1$ and $d_2$ are constants to be specified later. We focus on the Schnakenberg model (i.e. $g_1(u, v) = c(u - a - u^2v)$, $g_2(u, v) = c(u^2v-b)$) and the Gray-Scott model (i.e. $g_1(u, v) = (\mu + \rho)u - vu^2$, $g_2(u, v) = vu^2 - \rho(1-v)$), where $a, b, c, \mu$ and $\rho$ are constants.

\subsubsection{The Schnakenberg model} 
\begin{figure}[!ht]
\centering
\subfigure[$u(x)$]{\includegraphics[width=5.5cm]{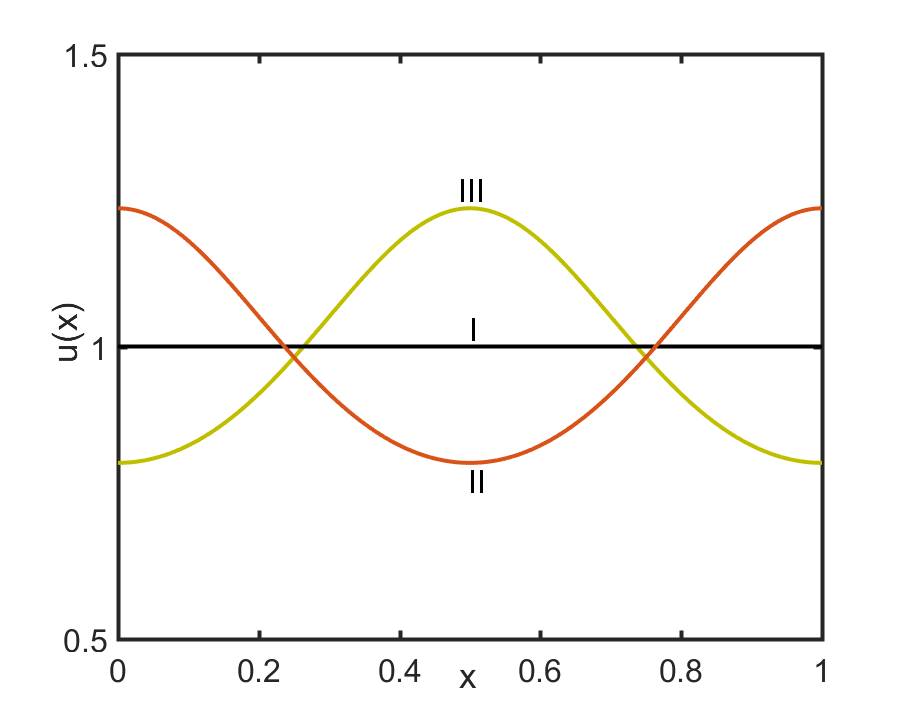}}\;\quad\quad
\subfigure[$v(x)$]{\includegraphics[width=5.5cm]{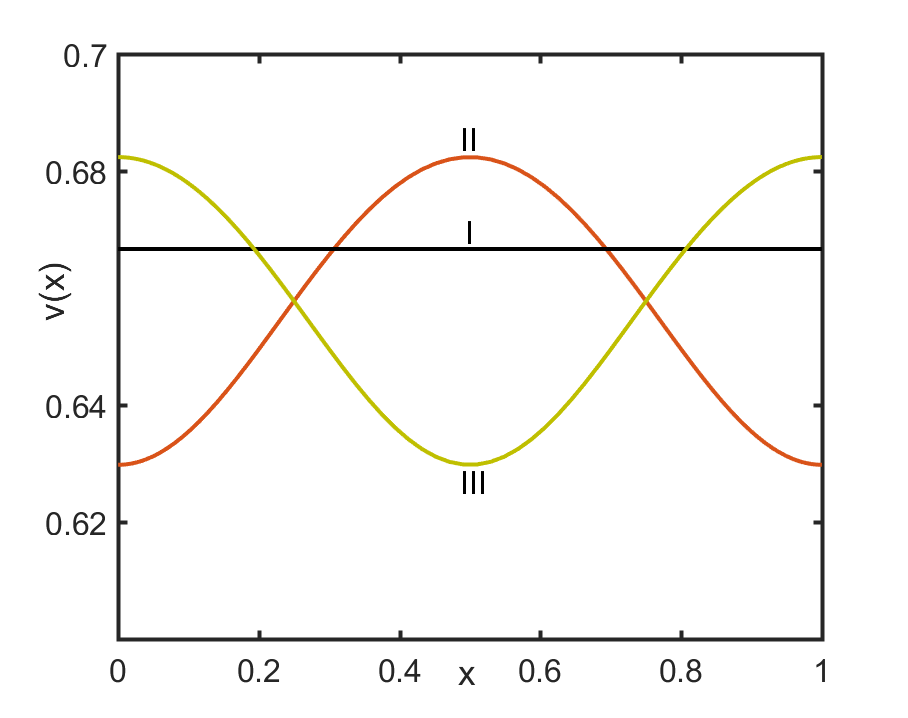}}\;
\caption{Multiple solutions of the Schnakenberg model with $d_2 = 50.$}\label{Figure3.3.1}
\end{figure}

\begin{figure}[!ht]
\centering
\subfigure[$u(x)$]{\includegraphics[width=5.5cm]{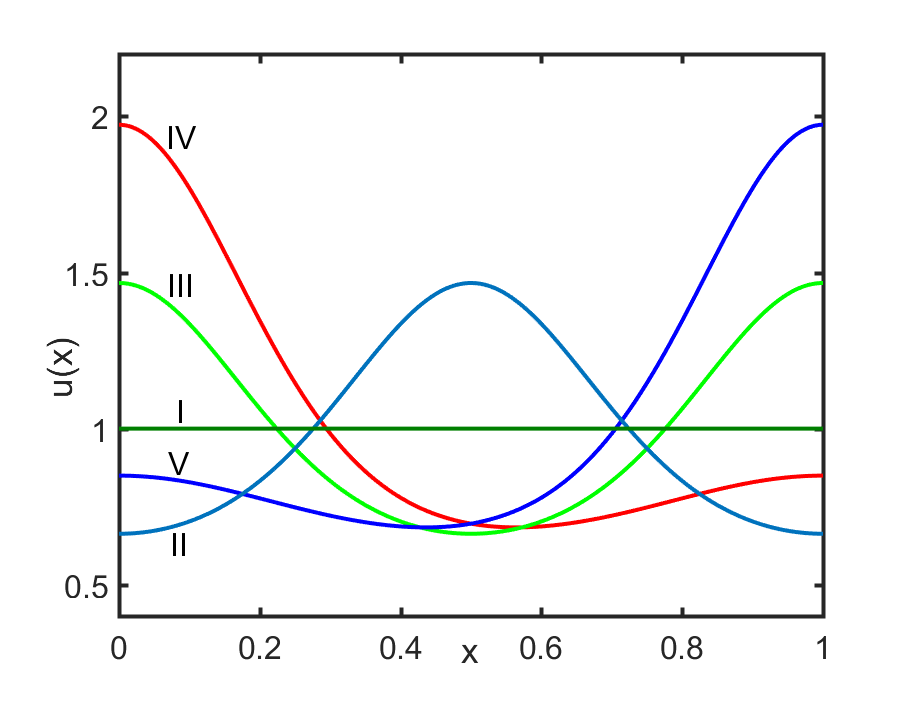}}\;\quad\quad
\subfigure[$v(x)$]{\includegraphics[width=5.5cm]{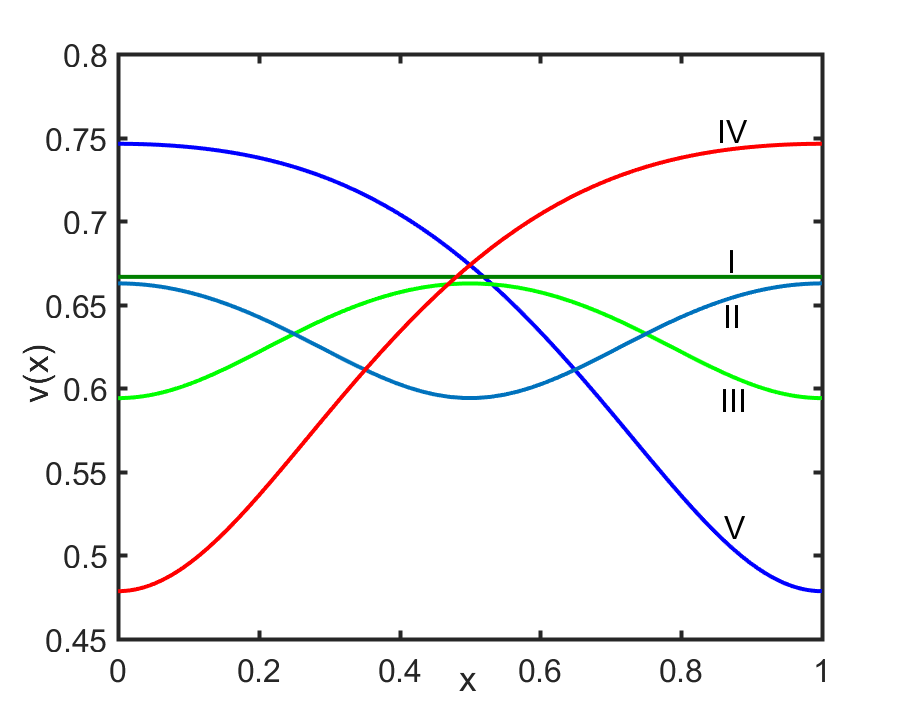}}\;
\caption{Multiple solutions of the Schnakenberg model with $d_2 = 70.$}\label{NewFigure3.3.1}
\end{figure}

As said in \cite{2020Spatial}, the Schnakenberg model is a Turing model, where $u$ is an activator and $v$ is a substrate. With numerical tests, we consider two cases $d_{2} =50$ and $d_{2} = 70$ with fixed parameters $d_1 = 1, a = 1/3, b = 2/3$ and $c = 200$, respectively. In Figs.\ref{Figure3.3.1}-\ref{NewFigure3.3.1}, multiple solutions of the Schnakenberg model on $u(x)$ and $v(x)$ are presented, where we list the number of solutions labeled by $\textrm{I, II}$, $\cdots$, $\textrm{V}$. These multiple solutions are agreement with that presented in \cite{2020Spatial}.
The efficiency of our spectral trust region LM-Deflation method for finding these multiple solutions is considered emphatically, and the corresponding numerical results are also included in Tables \ref{NewNTable41}-\ref{NewTable41}. Some observations and highlights are as follows.

\begin{table}[!h]
\centering
\caption{\small \normalsize Accuracy of our spectral trust region LM-Deflation method.}
\label{NewNTable41}\small
\setlength{\tabcolsep}{1.1mm}{
\begin{tabular}{c c c c c|| c c c c c c}
	\hline
~   &\multicolumn{4}{c}{$d_{2} = 50$} &\multicolumn{6}{c}{$d_{2} = 70$}\\\cline{2-11}
~  &$N$ &I  &II   &III  &$N$ &I  &II   &III  &IV   &V \\\cline{2-11}
\multirow{3}{*}{$||u_{N} - \hat{u}||_{\infty}$} &8   &1.06e-14   &5.78e-2      &2.06e-2     & 8  & 3.03e-13  & 1.97e-3   & 3.92e-2  &2.57e-3   &1.49e-3\\
~ &16  &2.94e-15   &3.29e-7      &9.20e-6    & 16  & 7.96e-15  & 5.82e-8   & 8.37e-6  &4.58e-8   &7.47e-8\\
~ &24  &2.37e-15   &8.14e-13     &5.21e-12  & 24   & 6.35e-15  & 9.16e-13  & 1.46e-13 &6.79e-13   &6.73e-14\\\hline\hline
\multirow{3}{*}{$||v_{N} - \hat{v}||_{\infty}$} &8 &4.69e-13 &3.96e-3 &8.34e-2 &8  &5.82e-13  &5.89e-3 &2.71e-2 &2.07e-3 &9.45e-3\\
~ &16  &7.93e-14   &5.20e-7      &4.62e-5    & 16  &9.26e-15  & 6.92e-8   &4.57e-6  &2.69e-9   &8.41e-8\\
~ &24  &1.02e-15   &6.74e-13     &3.02e-12  & 24   & 4.15e-15  & 8.31e-15  &8.47e-14 &9.36e-12   &3.24e-13\\\hline
\end{tabular}}
\end{table}

\begin{table}[!h]
\centering\small
\caption{\small A comparison of $\|\textbf{\textit{F}}(\textbf{\textit{x}})\|_{2}$ between our method and other methods.}
\label{NewTable41}
\setlength{\tabcolsep}{1.1mm}{
\begin{tabular}{cccc|cccc|cccc}
	\hline
\multicolumn{4}{c}{Newtonian iteration in \cite{farrell2015deflation}} &\multicolumn{4}{c}{LSTR method in \cite{2022Two}} &\multicolumn{4}{c}{our method}\\\hline
$n_{it}$ & $IG_{1}$  & $IG_{2}$   & $IG_{3}$   &$n_{it}$ & $IG_{1}$  & $IG_{2}$   & $IG_{3}$   &$n_{it}$ & $IG_{1}$  & $IG_{2}$   & $IG_{3}$   \\
5  &2.78e14  &3.81e13   &4.98e16    &5  &5.03e9  &3.20e7   &4.08e10  &5  &2.08e8  &5.21e10   &9.38e6 \\
10  &4.51e17 &9.06e18   &5.92e19   &10  &4.78e5 &3.01e3   &5.78e7   &10  &2.51e5 &9.38e5    &4.78e-1 \\
15  &3.01e23 &9.43e23   &8.49e26  &15  &3.92e1 &8.29e-1  &6.93e2   &15  &7.53e-1 &2.04e-2  &4.80e-7 \\
20  &4.83e32 &9.50e31   &8.71e30 &20  &3.19e-5 &9.85e-7 &4.32e-5  &20  &5.97e-6 &3.01e-7  &5.83e-10 \\
~  &- &-   &-                    &25  &5.29e-10 &6.05e-10 &3.50e-10  &25  &6.35e-11 &3.01e-13  &5.49e-13 \\\hline
T &- &-   &-               &~  &4.82 &3.72 &5.73  &~  &3.02 &2.93  &3.01 \\\hline
\end{tabular}}
\end{table}

\begin{itemize}
  \item When $N = 8, 16$ and 24, the accuracies of our spectral trust region LM-Deflation method are presented in Table \ref{NewNTable41}, where we compare the errors in maximum norm with the numerical solution $(\hat{u}, \hat{v})$ obtained with a relatively large $N$, indicating that our method is quite accurate.
  \item Based on the following initial guesses, a comparison of $\|\textbf{\textit{F}}(\textbf{\textit{x}})\|_{2}$ between our method and other methods is given in Table \ref{NewTable41}.
      \begin{equation}
      IG_{1}:
      \begin{cases}
      \tilde{u}^{(0)} = -\textrm{ones}(N+1, 1)\\
      \tilde{v}^{(0)} = -\textrm{ones}(N+1, 1)
      \end{cases}
      \;
      IG_{2}:
      \begin{cases}
      \tilde{u}^{(0)} = -\sin(\textrm{ones}(N+1, 1))\\
      \tilde{v}^{(0)} = -\textrm{ones}(N+1, 1)
      \end{cases}
      \;
      IG_{3}:= \sin(IG_{1}),
      \end{equation}
      where $N = 24$. The Newtonian iterations fail to converge for these initial inputs, while for LSTR method and our method, the values of $\|\textbf{\textit{F}}(\textbf{\textit{x}})\|_{2}$ descends very fast. However, it is noteworthy that the computational times (denoted by T) of our method are much less than that of LSTR method, indicating that the trust region Levenberg-Marquardt iteration in section \ref{sect2.3} plays an important role in improving the efficiency of our method.
\end{itemize}

\subsubsection{The Gray-Scott model} In the chemical reaction, the Gray¨CScott model is used to describe the autocatalytic reaction between the activator and the substrate. With the bootstrapping method, Hao et al. \cite{2020Spatial} studied multiple solutions of stationary spatial patterns on the Gray-Scott model on a 1D domain $[0, 1]$. Our spectral trust region LM-Deflation method is also used to find these multiple solutions, and we choose $d_1 = 2.5\times 10^{-4}, d_2 = 5\times 10^{-4}, \rho = 0.04$ and $\mu = 0.065$, which is the same as in \cite{2020Spatial}. Moreover, the property of the solutions on the Gray-Scott model is introduced, i.e. the solutions of the Gray-Scott model are symmetric with respect to the center of the domain $x = 0.5$. To be specific, if $u(x)$ is a solution, then $u((x+0.5)\; \textrm{mod}\; 1)$ is also a solution.

\begin{figure}[!ht]
\centering
\subfigure[I]{\includegraphics[width=2.5cm]{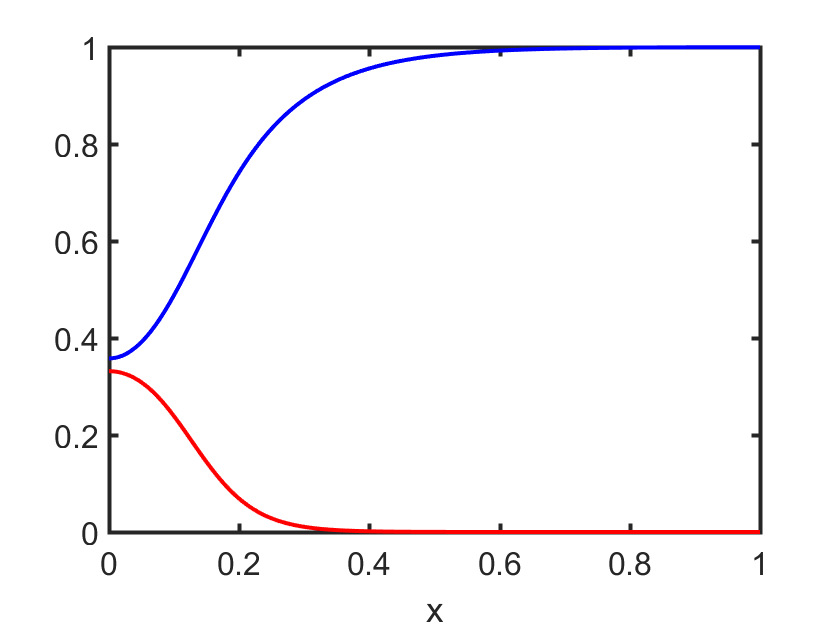}}\;
\subfigure[II]{\includegraphics[width=2.5cm]{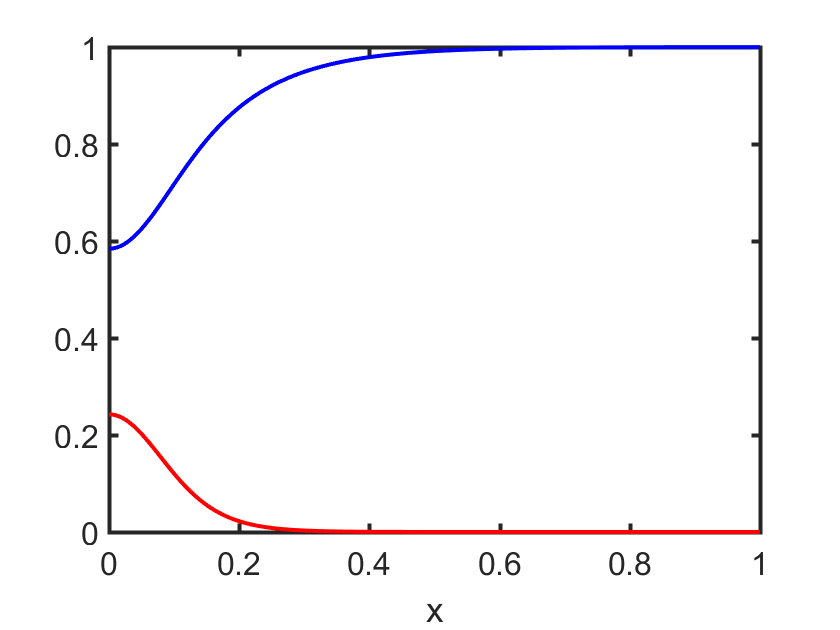}}\;
\subfigure[III]{\includegraphics[width=2.5cm]{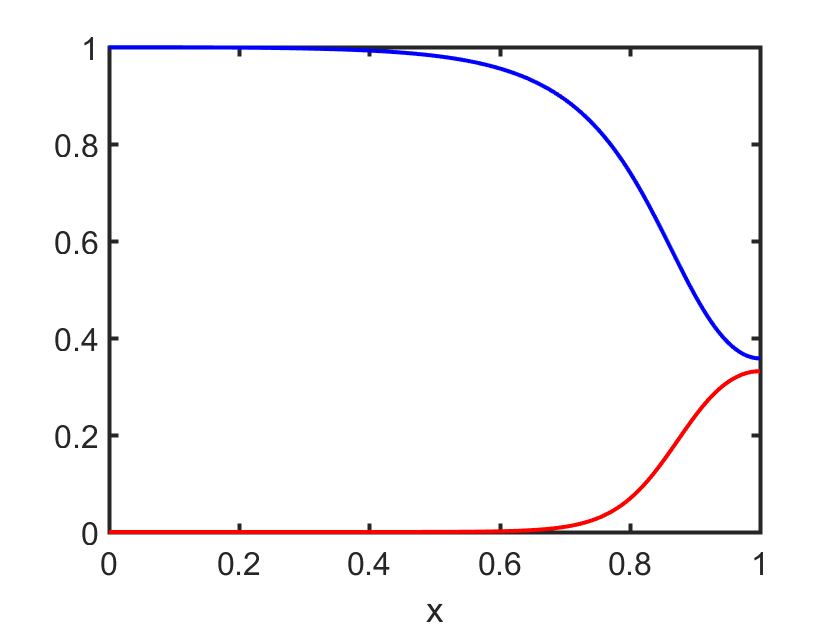}}\;
\subfigure[IV]{\includegraphics[width=2.5cm]{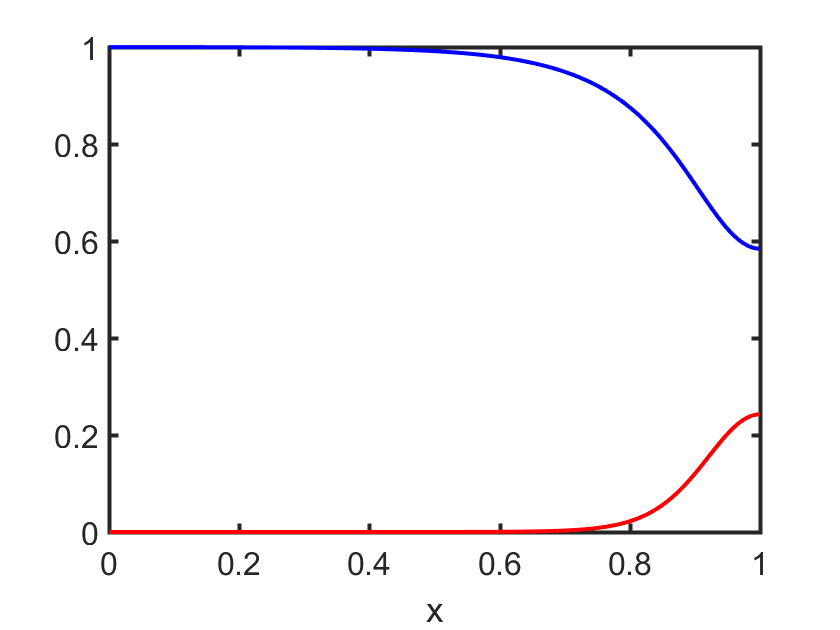}}\;\\\vspace{-0.3cm}
\subfigure[V]{\includegraphics[width=2.5cm]{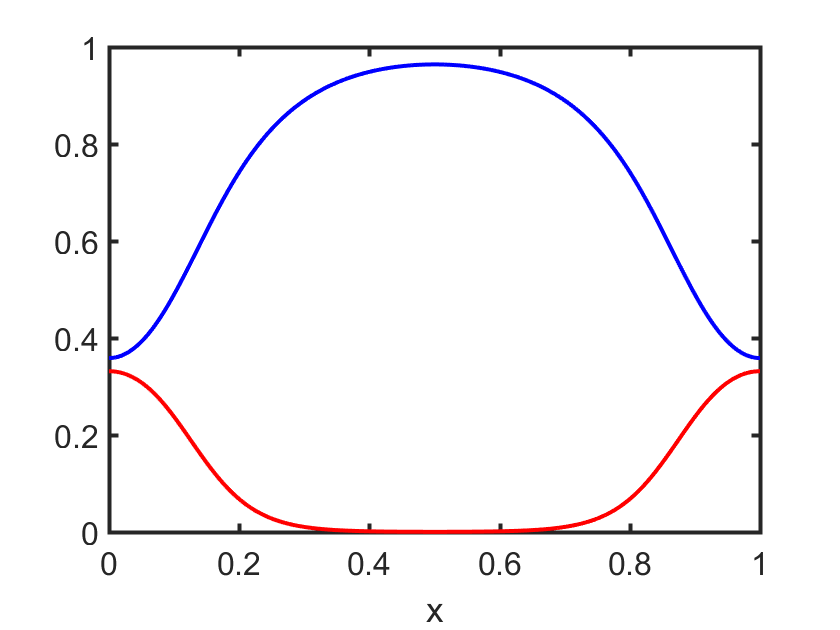}}\;
\subfigure[VI]{\includegraphics[width=2.4cm]{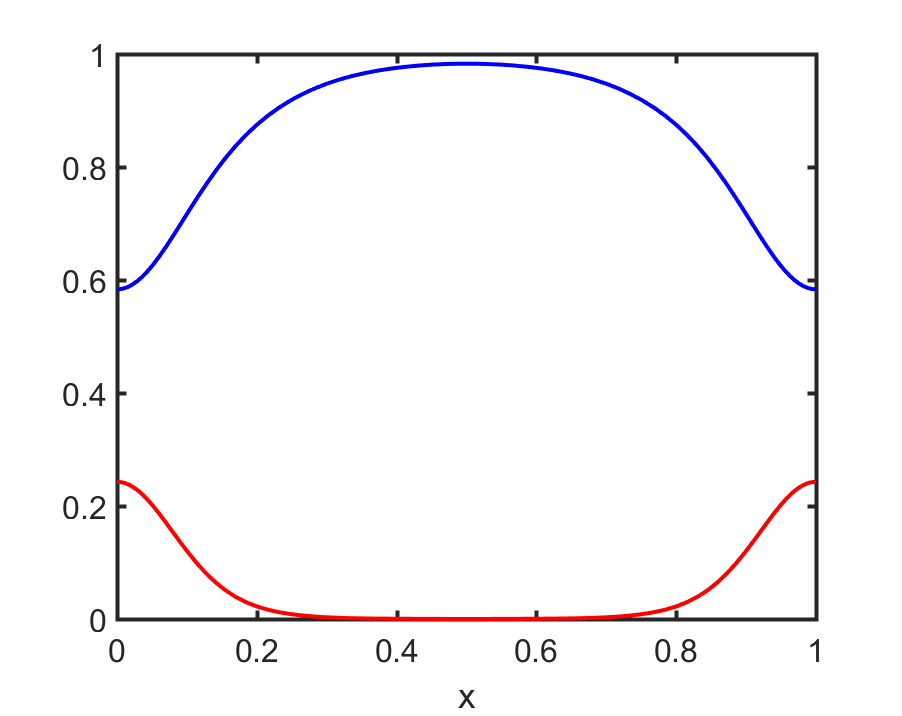}}\;
\subfigure[VII]{\includegraphics[width=2.4cm]{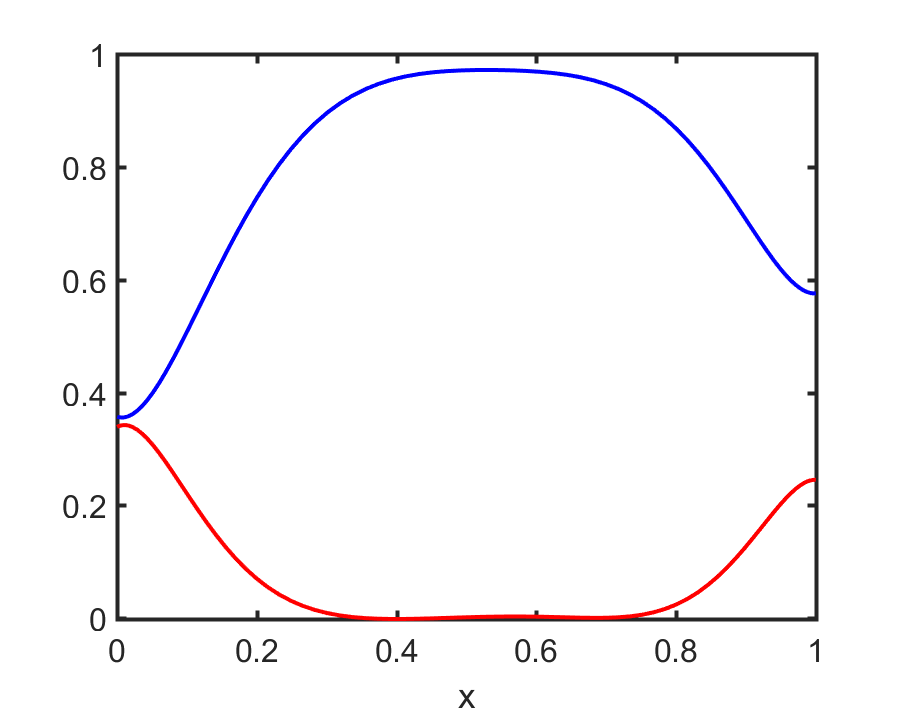}}\;
\subfigure[VIII]{\includegraphics[width=2.5cm]{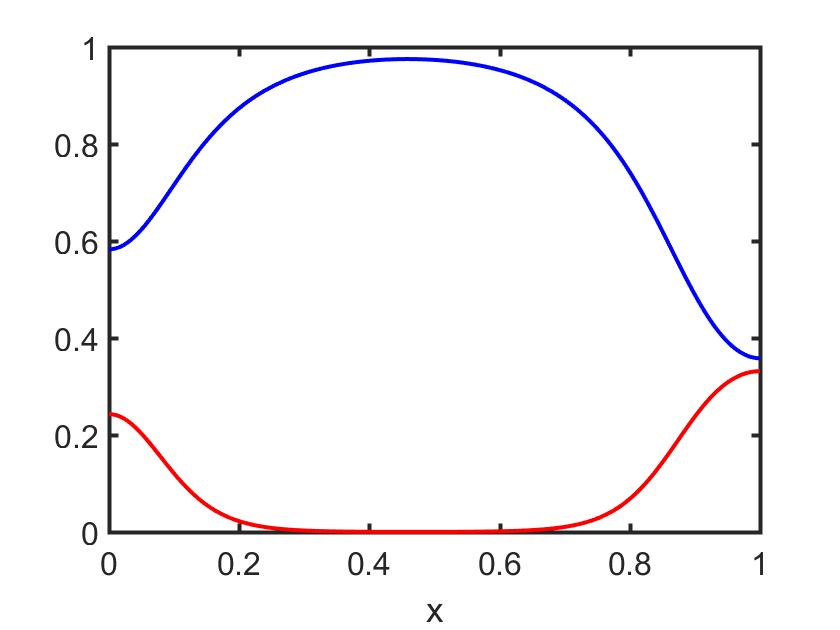}}\;\\\vspace{-0.3cm}
\subfigure[IX]{\includegraphics[width=2.5cm]{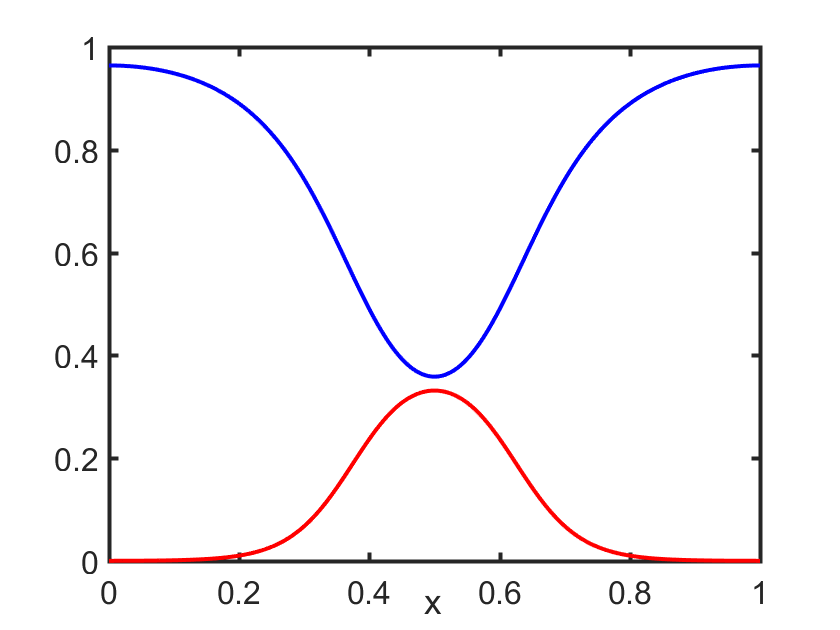}}\;
\subfigure[X]{\includegraphics[width=2.5cm]{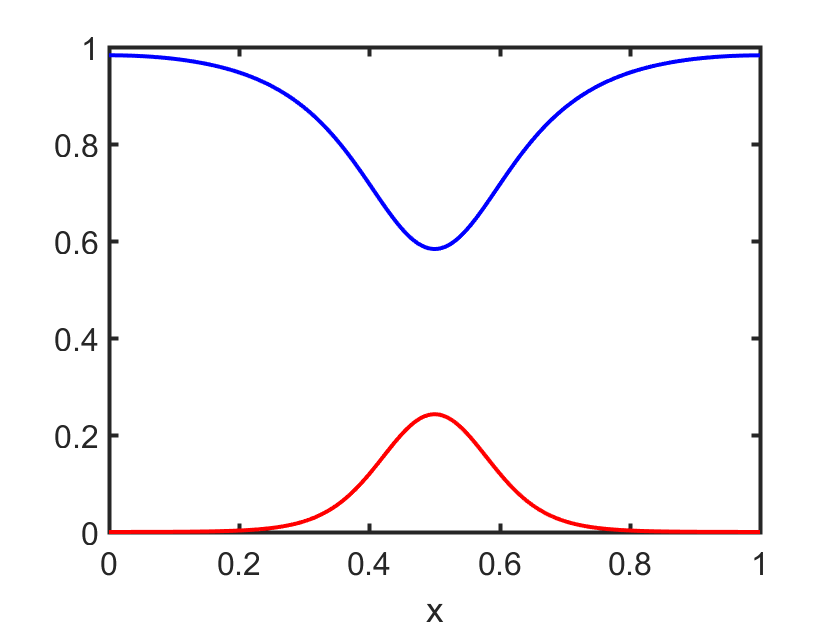}}\;
\subfigure[XI]{\includegraphics[width=2.5cm]{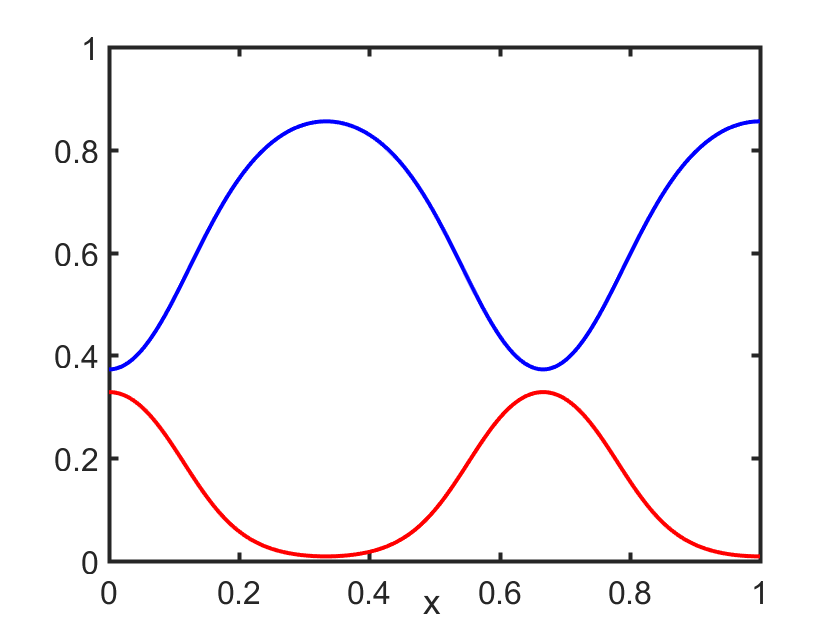}}\;
\subfigure[XII]{\includegraphics[width=2.5cm]{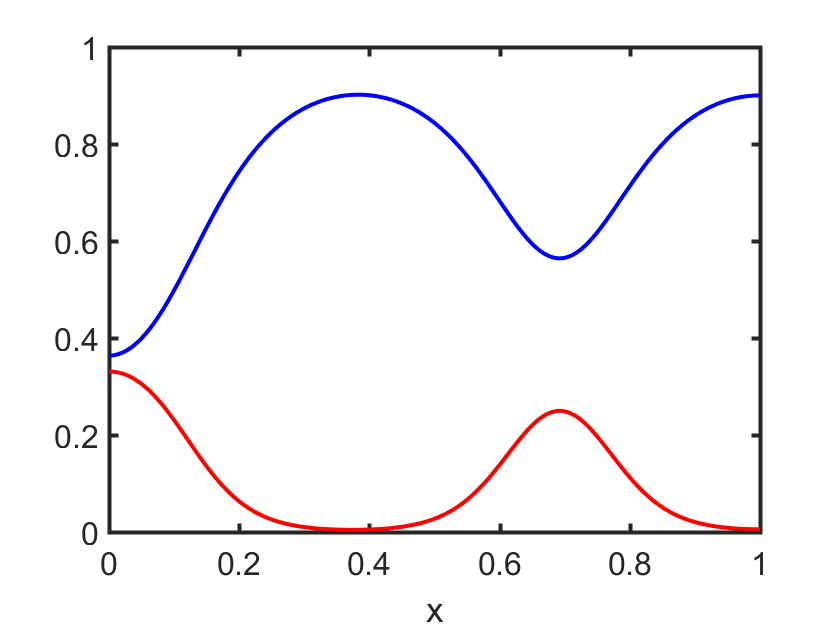}}\;\\\vspace{-0.3cm}
\subfigure[XIII]{\includegraphics[width=2.5cm]{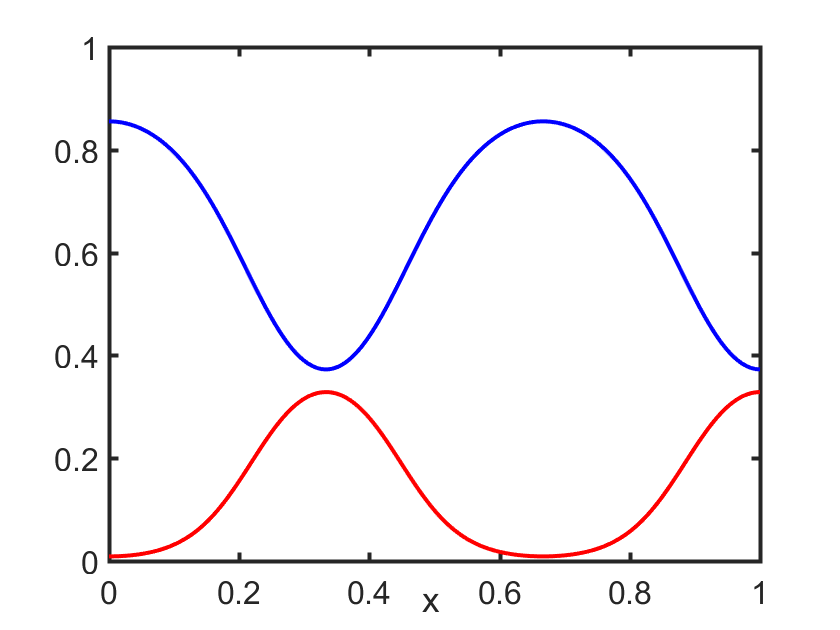}}\;
\subfigure[XIV]{\includegraphics[width=2.5cm]{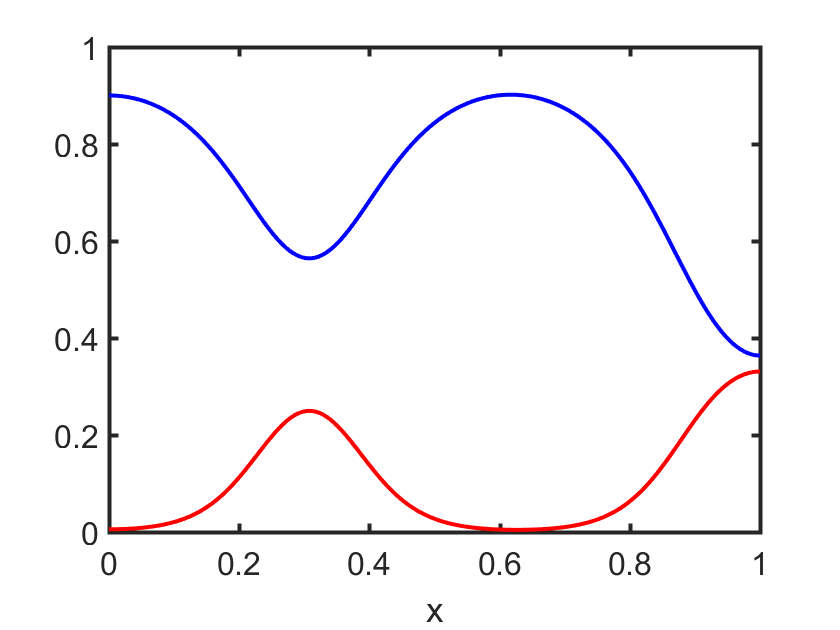}}\;
\subfigure[XV]{\includegraphics[width=2.5cm]{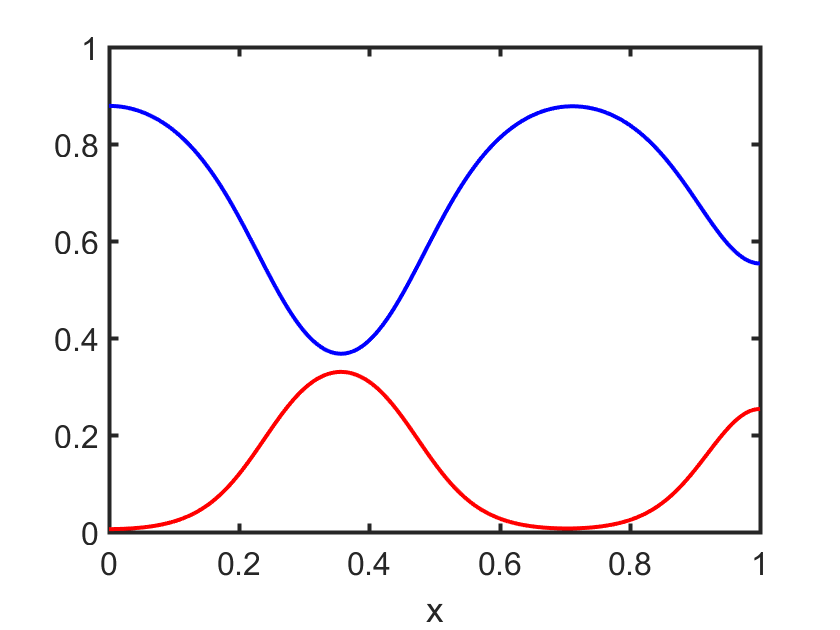}}\;
\subfigure[XVI]{\includegraphics[width=2.5cm]{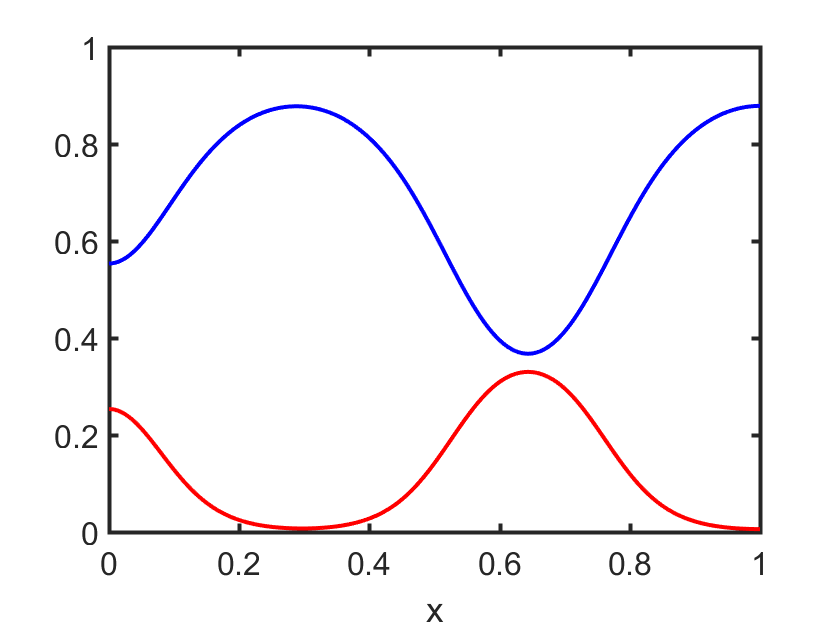}}\;\\\vspace{-0.3cm}
\subfigure[XVII]{\includegraphics[width=2.5cm]{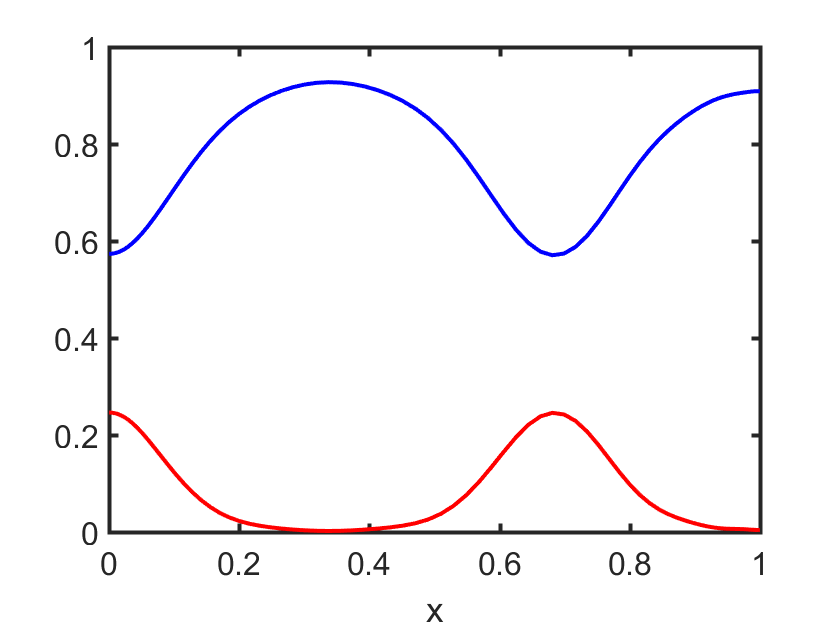}}\;
\subfigure[XVIII]{\includegraphics[width=2.5cm]{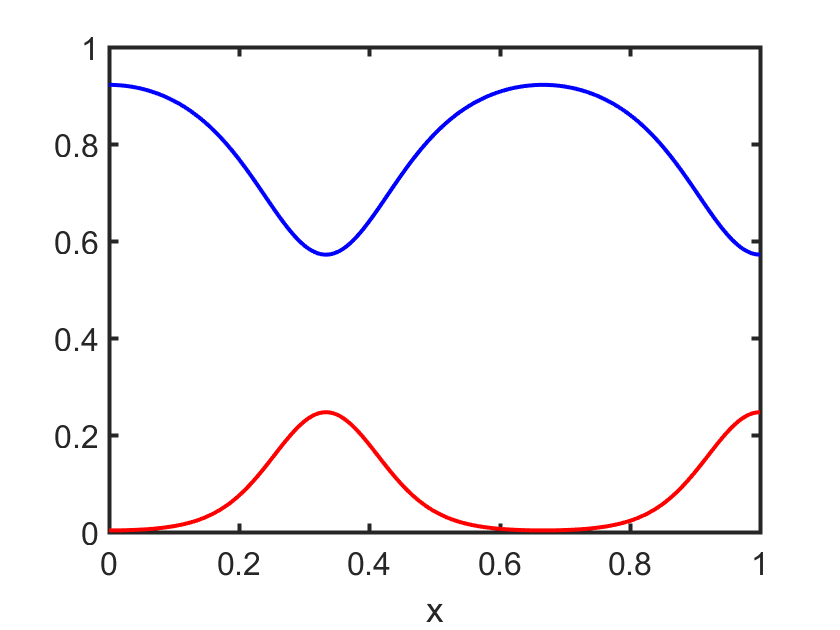}}\;
\subfigure[XIX]{\includegraphics[width=2.5cm]{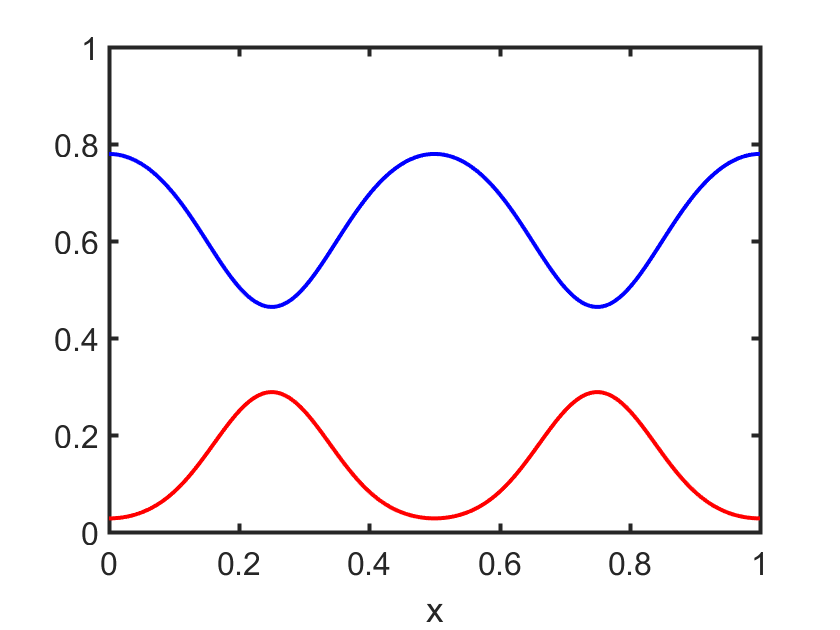}}\;
\subfigure[XX]{\includegraphics[width=2.5cm]{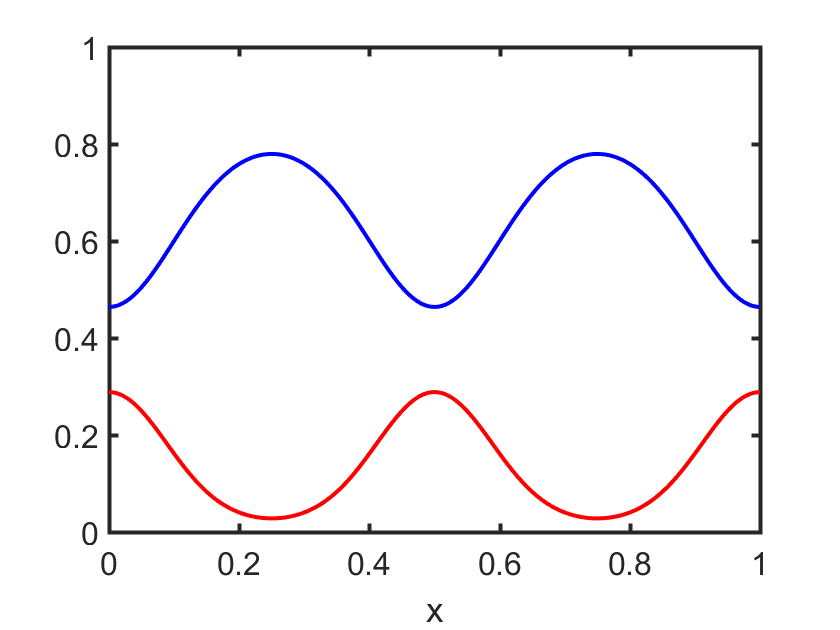}}\;
\vspace{-0.3cm}
\caption{Multiple solutions of the Gray-Scott model (Red: $u_{24}(x)$; Blue: $v_{24}(x)$).}\label{Example2.1}
\end{figure}

\begin{table}[!h]
\centering
\caption{\small \normalsize Accuracy of our spectral trust region LM-Deflation method on $||u_{N} - \hat{u}||_{\infty}$.}
\label{modifytable21}\small
\setlength{\tabcolsep}{1.1mm}{
\begin{tabular}{c c c c c c c c c c c c}
	\hline
 $N$ &I  &II   &III    &IV   &V  &VI   &VII  &VIII  &IX   &X \\\hline
8  &1.03e-3   & 5.21e-2   & 5.93e-2   &6.30e-2   & 6.29e-4   &3.92e-3   & 6.55e-2   &1.71e-3   &3.18e-2   &7.06e-2\\
16 &7.57e-6   &4.45e-7   &4.89e-7   &6.46e-6   &9.50e-8   &3.44e-8   &1.86e-7   &3.44e-6   &7.91e-8   &3.17e-7\\
24 &5.73e-10   &4.38e-11   &6.93e-12   &5.93e-11   &9.52e-11   &1.03e-12   &8.41e-10   &6.55e-11   &9.21e-13   & 2.09e-12\\\toprule[0.8pt]
$N$ &XI  &XII   &XIII    &XIV   &XV  &XVI   &XVII   &XVIII    &XIX    &XX \\\hline
8  &2.76e-2   & 4.98e-2   &9.68e-3   &8.42e-2   &6.99e-3   &1.63e-3   &9.59e-2   &6.55e-3   &1.19e-2   &2.43e-3\\
16 &5.43e-5   &6.83e-6   &2.54e-7   &8.14e-7   &9.29e-8   &2.55e-7   &5.47e-7   &8.14e-6   &8.28e-7   &7.34e-6\\
24 &5.83e-13   &1.65e-12   &2.67e-10   &6.66e-13   &4.98e-12   &1.38e-11   &9.59e-12   &7.51e-13   &2.53e-12   &2.08e-12\\\hline
\end{tabular}}
\end{table}

\begin{table}[!h]
\centering
\caption{\small \normalsize Accuracy of our spectral trust region LM-Deflation method on $||v_{N} - \hat{v}||_{\infty}$.}
\label{modifytable22}\small
\setlength{\tabcolsep}{1.1mm}{
\begin{tabular}{c c c c c c c c c c c c}
	\hline
 $N$ &I  &II   &III    &IV   &V  &VI   &VII  &VIII  &IX   &X \\\hline
8  &3.49e-3   &3.51e-2   &7.58e-2   &1.29e-3   &1.96e-2   &8.30e-3   & 7.57e-2   &5.39e-3   &5.68e-3   &2.51e-3\\
16 &2.51e-8   &3.81e-7   &5.30e-8   &4.69e-6   &1.19e-8   &2.51e-7   &6.16e-8   &4.73e-7   &9.17e-7   &5.33e-6\\
24 &5.85e-13   &3.37e-11   &1.62e-12   &6.01e-13   &8.25e-12   &3.11e-11   &6.89e-13   &4.42e-11   &8.17e-12   & 7.74e-13\\\toprule[0.8pt]
$N$ &XI  &XII   &XIII    &XIV   &XV  &XVI   &XVII   &XVIII    &XIX    &XX \\\hline
8  &8.68e-3   & 4.31e-2   &1.36e-3   &8.53e-2   &7.59e-2   &8.44e-1   &8.69e-3   &3.50e-2   &2.39e-2   &1.23e-3\\
16 &5.77e-5   &5.13e-7   &3.50e-6   &1.36e-8   &3.99e-6   &1.81e-7   &5.79e-7   &1.44e-6   &1.83e-8   &6.22e-6\\
24 &2.64e-13   &4.90e-12   &9.02e-12   &4.96e-12   &1.31e-13   &2.34e-11   &3.53e-11   &8.21e-12   &4.30e-12   &9.64e-11\\\hline
\end{tabular}}
\end{table}

\begin{table}[!h]
\centering\small
\caption{\small A comparison of $\|\textbf{\textit{F}}(\textbf{\textit{x}})\|_{2}$ between our method and other methods for the Gray-Scott model.}
\label{modifytable23}
\setlength{\tabcolsep}{1.1mm}{
\begin{tabular}{cccc|cccc|cccc}
	\hline
\multicolumn{4}{c}{Newtonian iteration in \cite{farrell2015deflation}} &\multicolumn{4}{c}{LSTR method in \cite{2022Two}} &\multicolumn{4}{c}{our method}\\\hline
$n_{it}$ & $IG_{1}$  & $IG_{2}$   & $IG_{3}$   &$n_{it}$ & $IG_{1}$  & $IG_{2}$   & $IG_{3}$   &$n_{it}$ & $IG_{1}$  & $IG_{2}$   & $IG_{3}$   \\
5  &1.68e9  &5.47e10   &3.06e10   &5   &5.03e6  &3.20e7   &4.08e8   &5   &4.09e7  &4.81e7   &4.96e8 \\
10 &1.88e14 &4.86e16   &5.10e14   &10  &5.56e3  &3.81e3   &2.75e4   &10  &2.51e4 &3.85e4    &7.78e4 \\
15 &3.21e17 &8.17e16   &5.10e15   &15  &3.85e0  &8.09e-1  &7.93e0   &15  &4.83e1 &2.67e-1  &4.29e-1 \\
20 &6.86e23 &4.50e25   &3.68e20   &20  &8.41e-5 &4.56e-7  &4.52e-5  &20  &3.39e-6 &8.83e-7  &4.21e-6 \\
~  &- &-   &-                    &25   &5.93e-10&1.98e-10 &3.67e-10 &25  &9.12e-11 &9.36e-13  &7.73e-13 \\\hline
T &- &-   &-               &~  &5.09 &5.31 &5.12  &~  &4.21 &3.39  &4.19 \\\hline
\end{tabular}}
\end{table}

\begin{figure}[!ht]
\begin{center}
\subfigure[$\textrm{XXI}$]{ \includegraphics[width=4.5cm]{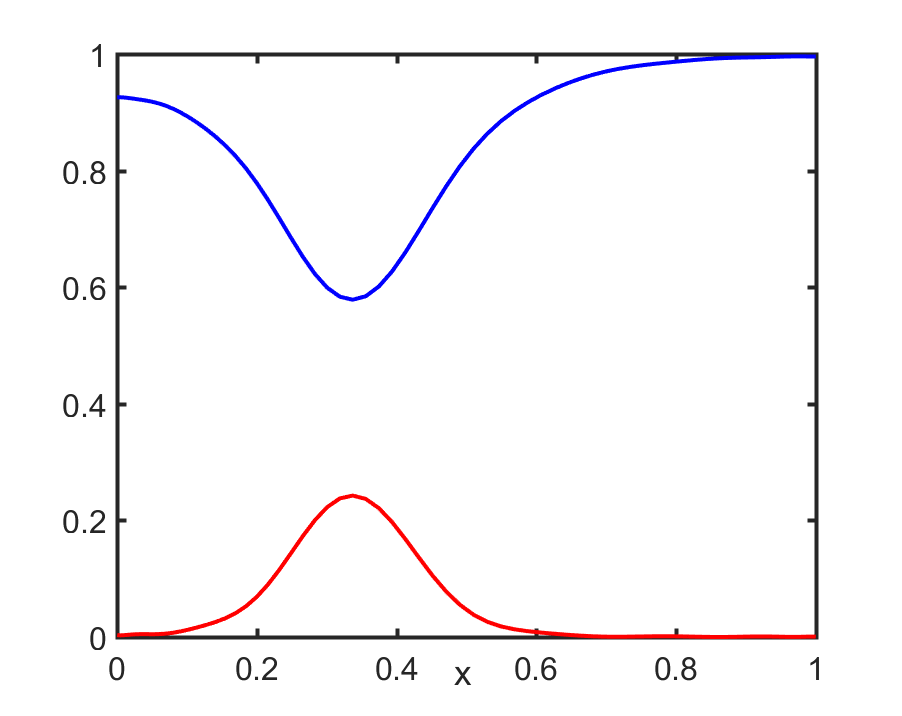}}\;\;
\subfigure[$\textrm{XXII}$]{ \includegraphics[width=4.5cm]{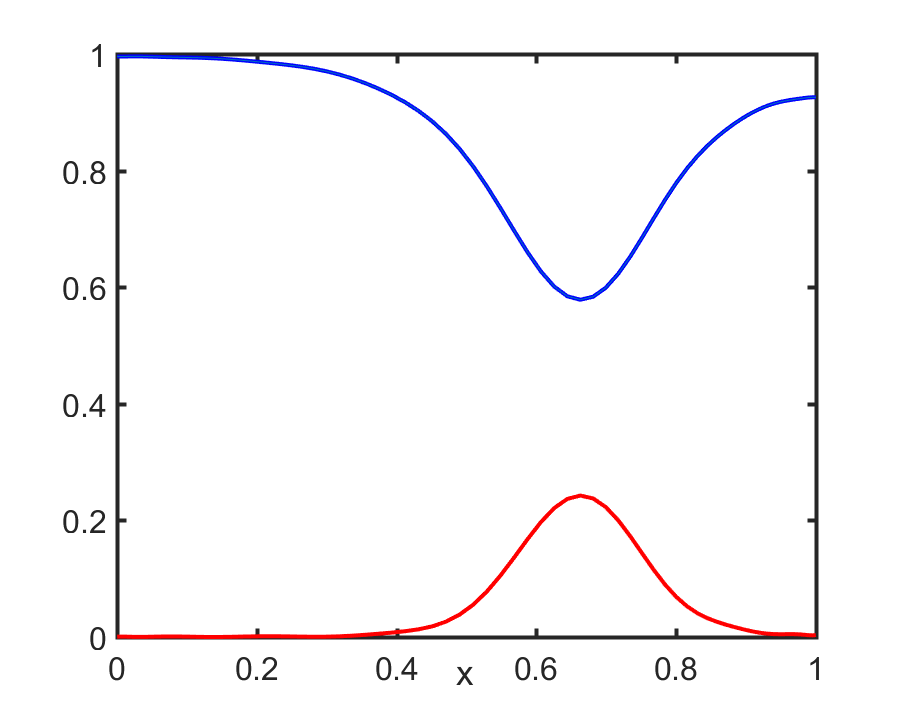}}\\\vspace{-0.4cm}
\caption{New solutions of the Gray-Scott model (Red: $u_{24}(x)$; Blue: $v_{24}(x)$).}
\label{Newexample1}
\end{center}
\end{figure}

In Fig.\ref{Example2.1}, we show multiple solutions of the Gray-Scott model in a 1D domain $[0, 1]$, and the number of solutions are labeled by $\textrm{I}, \textrm{II}, \cdots$, where these multiple solutions are agreement with that obtained in \cite{2020Spatial}, and the solutions are symmetric with respect to the center of the domain $x = 0.5$ (such as $\textrm{I}$ and $\textrm{III}$, $\textrm{II}$ and $\textrm{IV}$). Numerical tests on $\textrm{I}, \cdots, \textrm{XX}$ are considered to show the efficiency of our spectral trust region LM-Deflation method. In Tables \ref{modifytable21}-\ref{modifytable22}, with increasing $N$, the accuracies on $||u_{N} - \hat{u}||_{\infty}$ and $||v_{N} - \hat{v}||_{\infty}$ are given, indicating that our spectral trust region LM-Deflation method is reliable and effective. Since we are not trying to provide the Hessian matrix exactly, the efficiency of our spectral trust region LM-Deflation method in the nonlinear iteration is well behaved in Table \ref{modifytable23}. Furthermore, two new solutions which were not reported in the literature are also found (see Fig.\ref{Newexample1}), and they are labeled by $\textrm{XXI}$ and $\textrm{XXII}$. The corresponding numerical information is also given in Table \ref{modifytable24}.

\begin{table}[!h]
\centering\small
\caption{\small Performance of our spectral trust region LM-Deflation method to new solutions of the Gray-Scott model.}
\label{modifytable24}
\begin{tabular}{cccc|ccc}
  \hline
 ~ &\multicolumn{3}{c}{XXI} &\multicolumn{3}{c}{XXII}\\\cline{2-7}
  $N$           &8   &16     &24     &8      &16     &24 \\
  $L^{\infty}$ error &2.8401e-4 &3.9201e-8   &4.0923e-12   &3.6520e-5   &1.9302e-9   &1.2084e-13 \\
  $\|{\bs F}({\bs x})\|_{2}$ &4.0398e-12 &6.3028e-13   &6.9281e-12   &1.9026e-12   &7.9036e-12   &2.0538e-13 \\
  \hline
\end{tabular}
\end{table}

\subsection{PDE examples}\label{sect3.2} 
\begin{figure}[!ht]
\begin{center}
\subfigure[I]{ \includegraphics[width=5.5cm]{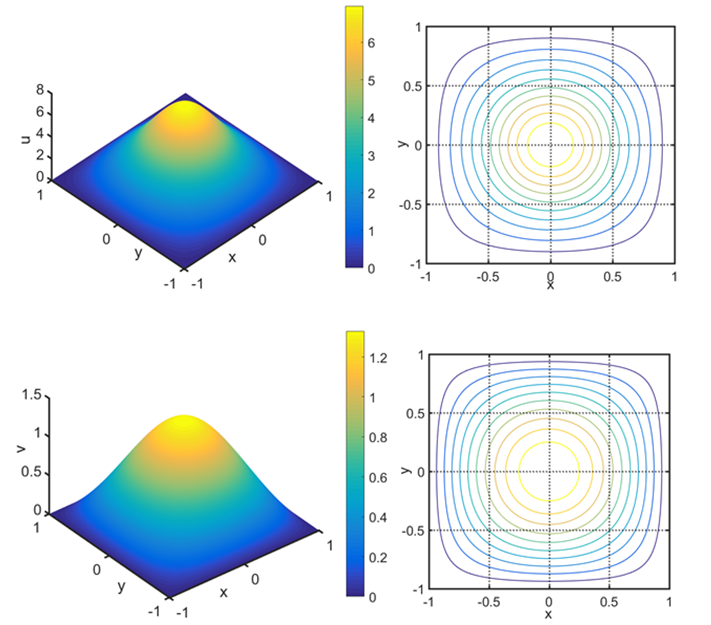}}\;
\subfigure[II]{ \includegraphics[width=5.5cm]{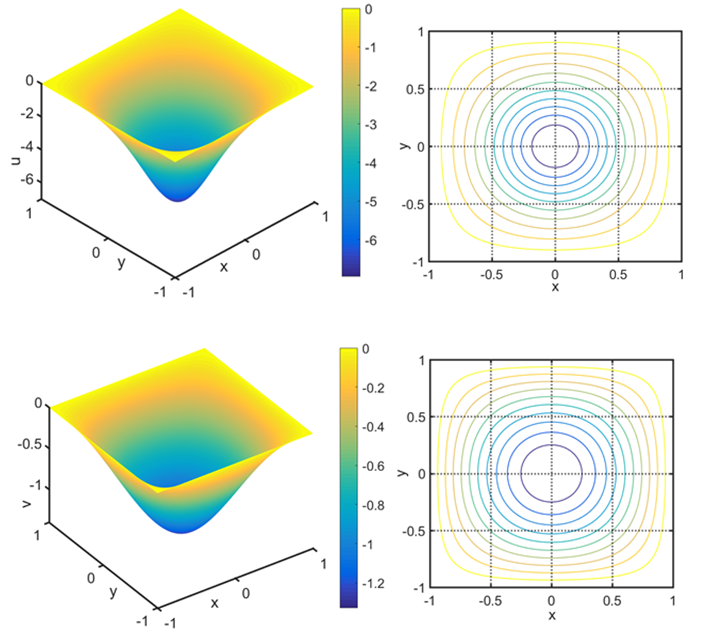}}\\
\caption{Profiles of 1-peak solutions for $p = 3$ and $q = 3$.}
\label{example3.1}
\end{center}
\end{figure}

\begin{figure}[!ht]
\begin{center}
\subfigure[III]{ \includegraphics[width=5.5cm]{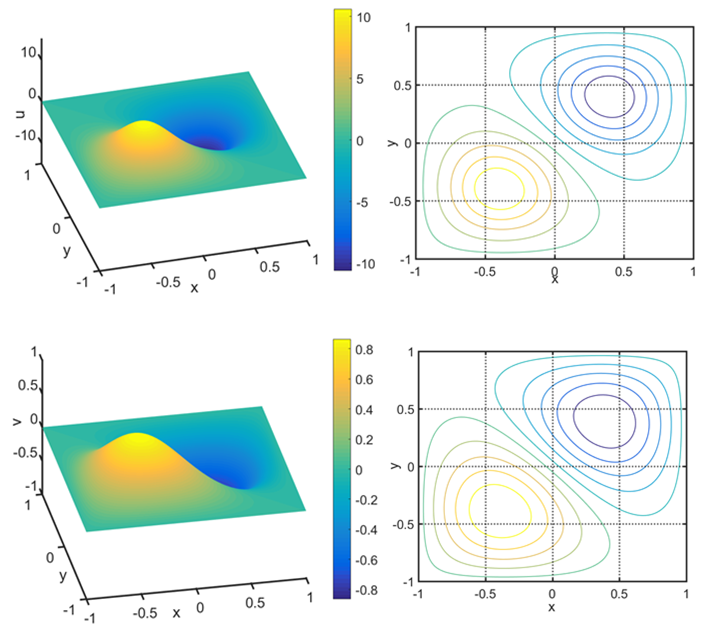}}\;
\subfigure[IV]{ \includegraphics[width=5.5cm]{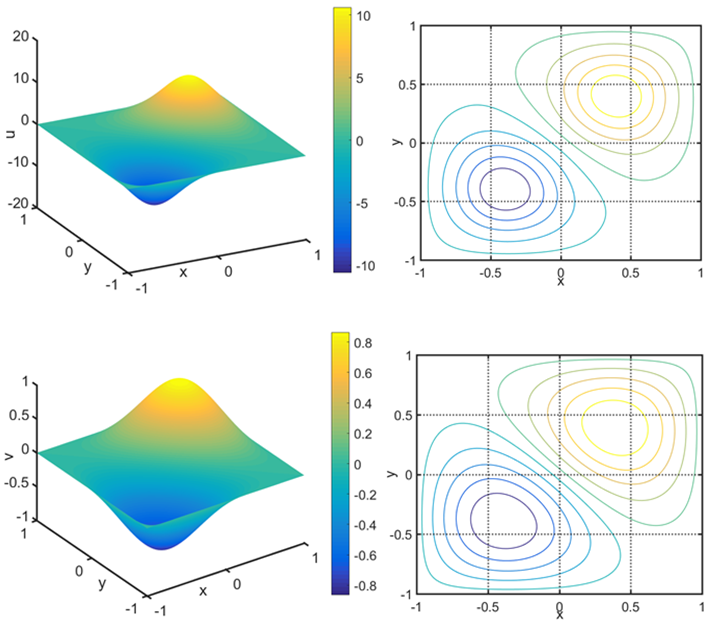}}\\
\caption{Profiles of 2-peak sign-changing solutions for $p = 3$ and $q = 3$.}
\label{example3.2}
\end{center}
\end{figure}

\begin{figure}[!h]
\begin{center}
\subfigure[V]{ \includegraphics[width=5.5cm]{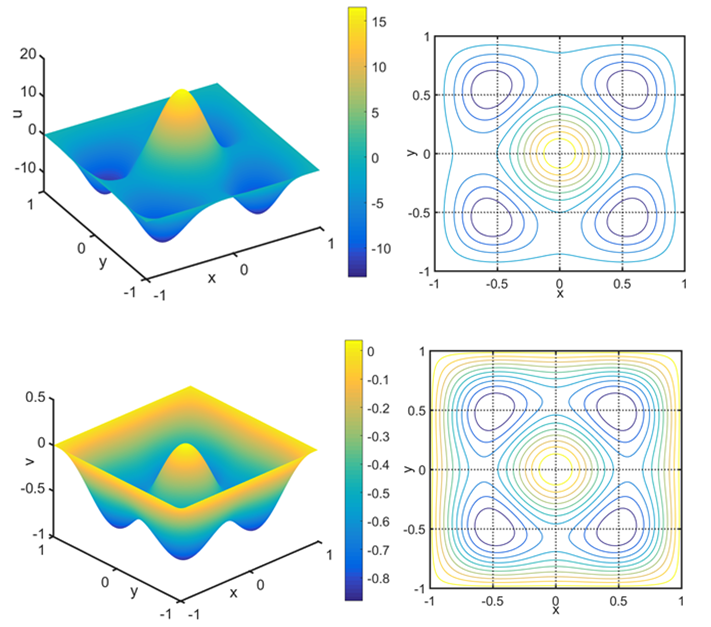}}\;
\subfigure[VI]{ \includegraphics[width=5.5cm]{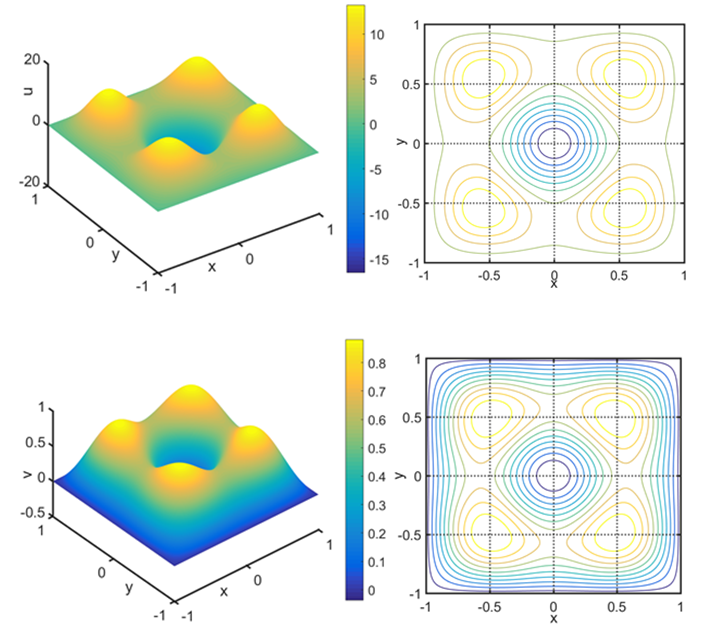}}\\
\caption{Profiles of multi-peak solutions for $p = 3$ and $q = 3$.}
\label{example3.3}
\end{center}
\end{figure}

We consider the model problem:
\begin{equation}\label{example3}
\begin{cases}
-\triangle u = G_1(u(x, y), v(x, y)),   \\
-\triangle v = G_2(u(x, y), v(x, y)),   \\
u = 0,\;   v = 0   \quad \textrm{on}\; \partial\Omega,
\end{cases}
\end{equation}
where $(x, y)\in \Omega$, and $G_1(u(x, y), v(x, y)), G_2(u(x, y), v(x, y))$ and $\Omega$ are to be specified later. Here we focus on three cases as follows:

\noindent\underline{{\bf Case 1:}} The problem \eqref{example3} with $G_1(u(x, y), v(x, y)) = \lambda u - \delta v + |u|^{p-1}u$ and $G_2(u(x, y), v(x, y)) = \delta u + \gamma v - |v|^{q-1}v$ is known as the noncooperative system of definite type \cite{2010A}, which together with the domain $\Omega = (-1, 1)\times(-1, 1)$. We choose $p = q = 3, \lambda = -0.5, \gamma = -0.5, \delta = 5.$ In Figs.\ref{example3.1}-\ref{example3.3}, we present multiple solutions of the noncooperative system of definite type and their contours, where we list the number of solutions labeled by $\textrm{I}, \textrm{II}, \cdots, \textrm{VI}$. The solutions with 1-peak or 2-peak are agreement with that presented in \cite{2010A}. While the solutions with multi-peak are found and shown for the first time. The accuracies of these multiple solutions are presented in Table \ref{PDE1modifytable1}. Based on different initial guesses in \eqref{PDE1} ($N = 24$), We further consider a comparison of $\|\textbf{\textit{F}}(\textbf{\textit{x}})\|_{2}$ between our method and other methods (see Table \ref{PDE1modifytable2}), indicating that our method is very effective again.

\begin{equation}\label{PDE1}
      IG_{1}:
      \begin{cases}
      \tilde{u}^{(0)} = -\textrm{ones}(N+1, N+1)\\
      \tilde{v}^{(0)} = -\textrm{ones}(N+1, N+1)
      \end{cases}
      \;
      IG_{2}:
      \begin{cases}
      \tilde{u}^{(0)} = -\sin(\textrm{ones}(N+1, N+1))\\
      \tilde{v}^{(0)} = -\sin(\textrm{ones}(N+1, N+1)),
      \end{cases}
      \end{equation}

\begin{table}[!h]
\centering
\caption{\small \normalsize Accuracy of our spectral trust region LM-Deflation method for $p = 3$ and $q = 3$.}
\label{PDE1modifytable1}\small
\setlength{\tabcolsep}{1.1mm}{
\begin{tabular}{c c c c c c c c}
	\hline
~   &$N$ &I  &II   &III    &IV   &V  &VI \\\cline{2-8}
\multirow{3}{*}{$||u_{24}-\hat{u}||_{\infty}$}&8 &6.92e-1 &7.65e-2 &4.97e-1 &1.62e-1 &2.55e-2 &7.09e-2\\
~  &16 &1.86e-5 &6.79e-6 &9.59e-7 &9.50e-6 &7.65e-6 &7.95e-5\\
~  &24 &4.38e-9 &3.81e-8 &5.85e-9 &5.05e-10 &4.45e-10 &1.62e-8\\\hline\hline
\multirow{3}{*}{$||v_{24}-\hat{v}||_{\infty}$}&8 &8.14e-1 &6.16e-2 &7.58e-1 &5.68e-2 &2.85e-2 &2.43e-2\\
~  &16 &1.96e-7 &5.67e-8 &2.51e-5 &9.34e-8 &1.19e-7 &4.69e-6\\
~  &24 &3.49e-10 &9.29e-11 &2.56e-10 &7.79e-9 &3.49e-8 &4.73e-9\\\hline
\end{tabular}}
\end{table}

\begin{table}[!h]
\centering\small
\caption{\small A comparison of $\|\textbf{\textit{F}}(\textbf{\textit{x}})\|_{2}$ between our method and other methods for $p = 3$ and $q = 3$.}
\label{PDE1modifytable2}
\setlength{\tabcolsep}{1.1mm}{
\begin{tabular}{ccc|ccc|ccc}
	\hline
\multicolumn{3}{c}{Newtonian iteration in \cite{farrell2015deflation}} &\multicolumn{3}{c}{LSTR method in \cite{2022Two}} &\multicolumn{3}{c}{our method}\\\hline
$n_{it}$ & $IG_{1}$  & $IG_{2}$   &$n_{it}$ & $IG_{1}$  & $IG_{2}$  &$n_{it}$ & $IG_{1}$  & $IG_{2}$\\
5 &2.24e10 &8.17e15          &5    &5.92e8  &8.32e10   &5   &2.51e9  &6.92e12 \\
10  &1.44e14  &8.38e19      &15   &2.91e5   &4.82e7    &15   &9.91e7  &6.93e8 \\
15  &9.13e20  &7.48e23      &25   &9.62e1   &7.71e2    &25   &4.73e3  &9.53e4 \\
20  &7.81e27  &7.74e29      &35   &3.19e-5  &8.50e-4   &35   &2.85e-3  &6.71e-4 \\
~  &- &-                    &45   &5.97e-10 &3.92e-10  &45   &4.52e-11 &7.83e-12 \\\hline
T &- &-                     &~    &9.34     &9.83      &~     &6.31    &5.38   \\\hline
\end{tabular}}
\end{table}

\begin{figure}[!ht]
\begin{center}
\subfigure[$(|u|u)'$ vs. $u$]{ \includegraphics[width=4.5cm]{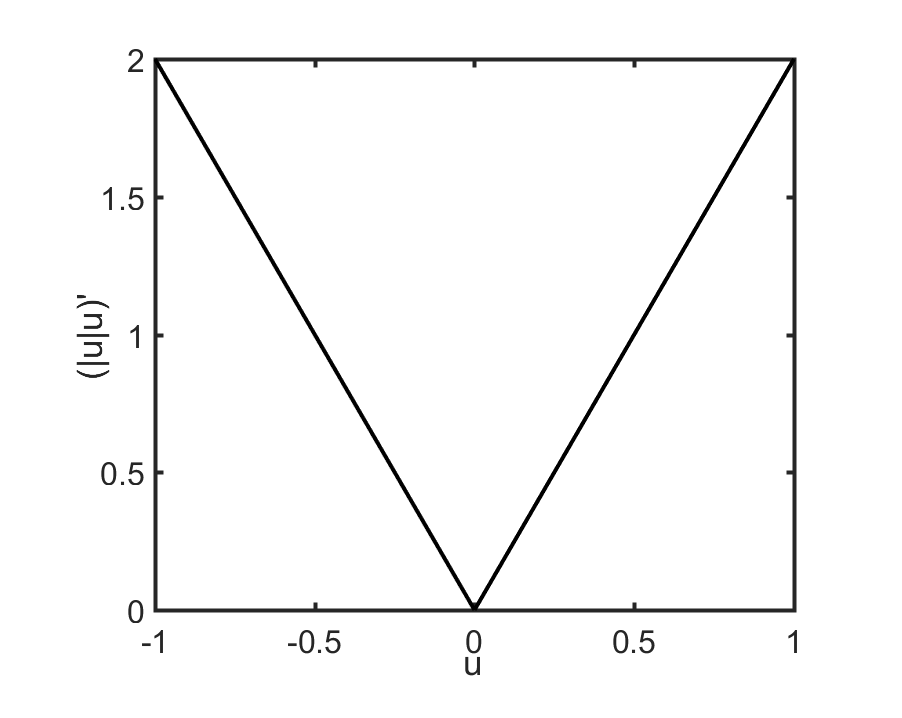}}\;
\subfigure[$(|u|u)''$ vs. $u$]{ \includegraphics[width=4.5cm]{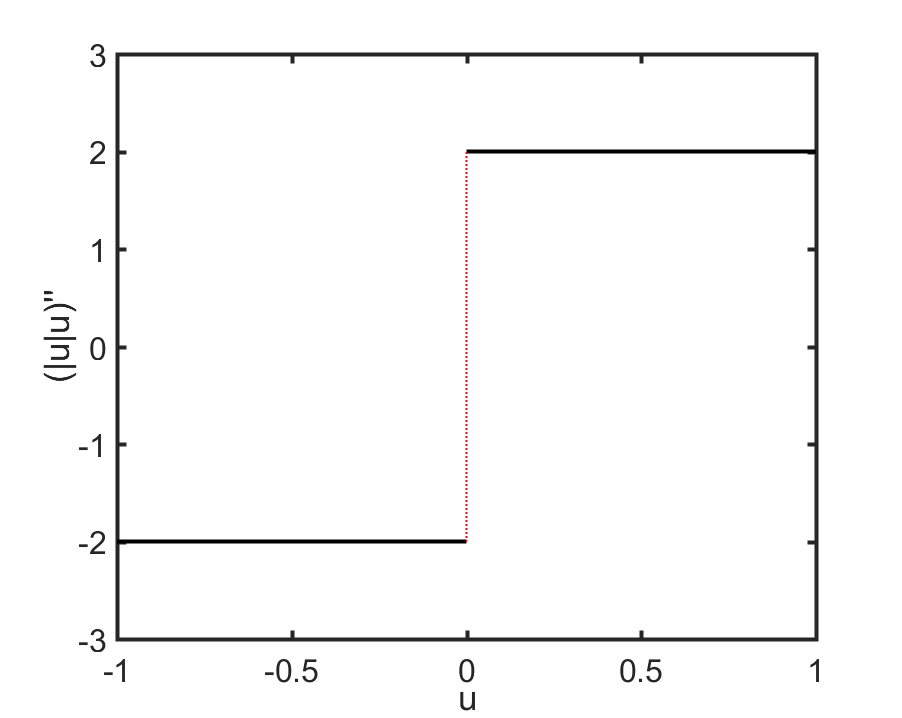}}\\
\caption{Profiles of $(|u|u)'$ and $(|u|u)''$ vs. $u$ on $[-1, 1]$.}
\label{pq2example1}
\end{center}
\end{figure}

\begin{figure}[!ht]
\begin{center}
\subfigure[I]{ \includegraphics[width=5.5cm]{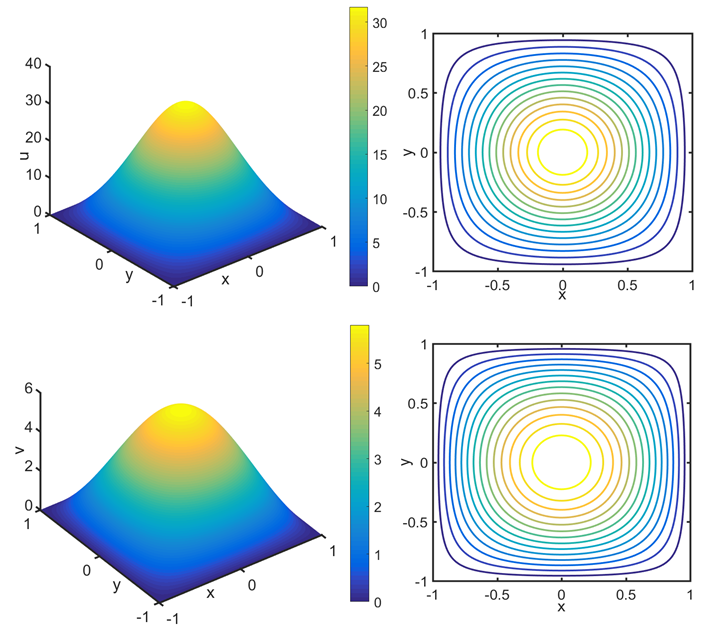}}\;
\subfigure[II]{ \includegraphics[width=5.5cm]{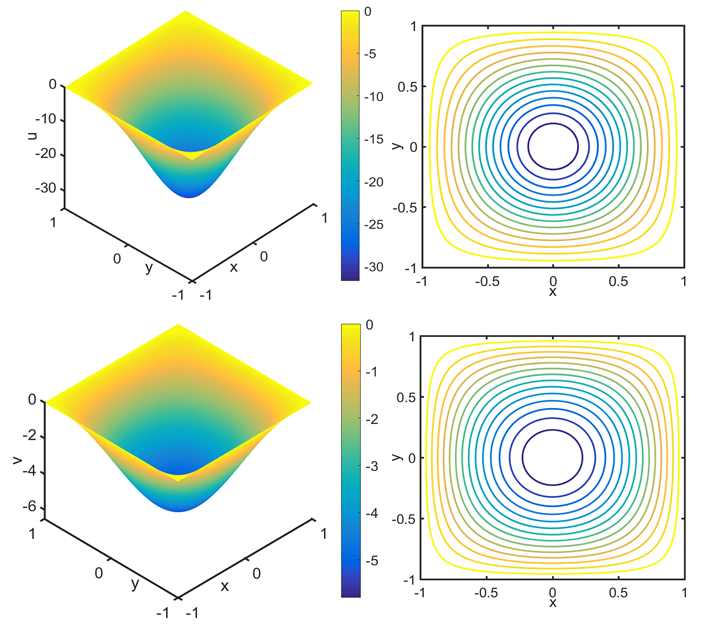}}\\
\caption{Profiles of 1-peak solutions for $p = 2, q = 2$.}
\label{pq2example2}
\end{center}
\end{figure}

\begin{figure}[!ht]
\begin{center}
\subfigure[III]{ \includegraphics[width=5.5cm]{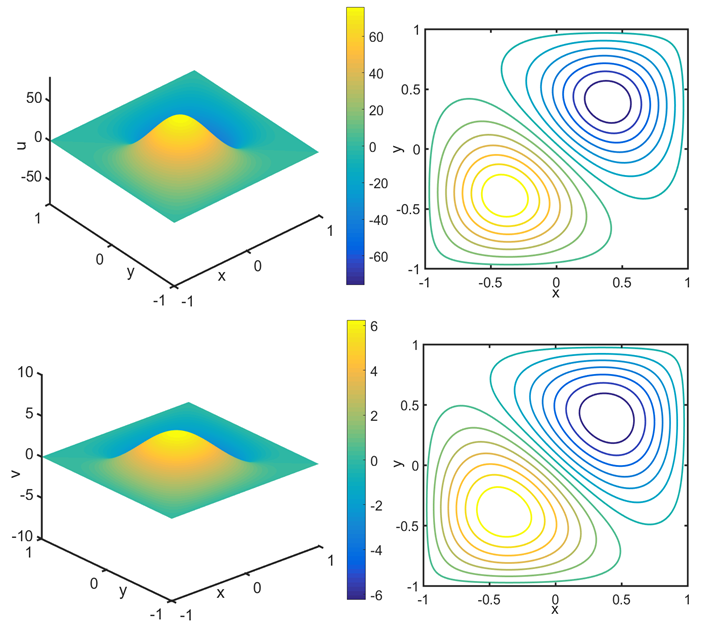}}\;
\subfigure[IV]{ \includegraphics[width=5.5cm]{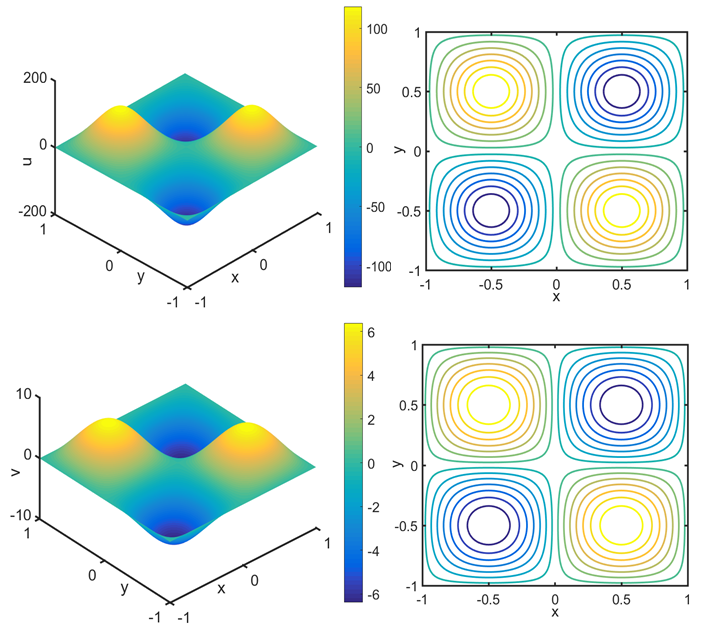}}\\
\caption{Profiles of 2-peak sign-changing solutions for $p = 2, q = 2$.}
\label{pq2example3}
\end{center}
\end{figure}

\begin{figure}[!ht]
\begin{center}
\subfigure[V]{ \includegraphics[width=5.5cm]{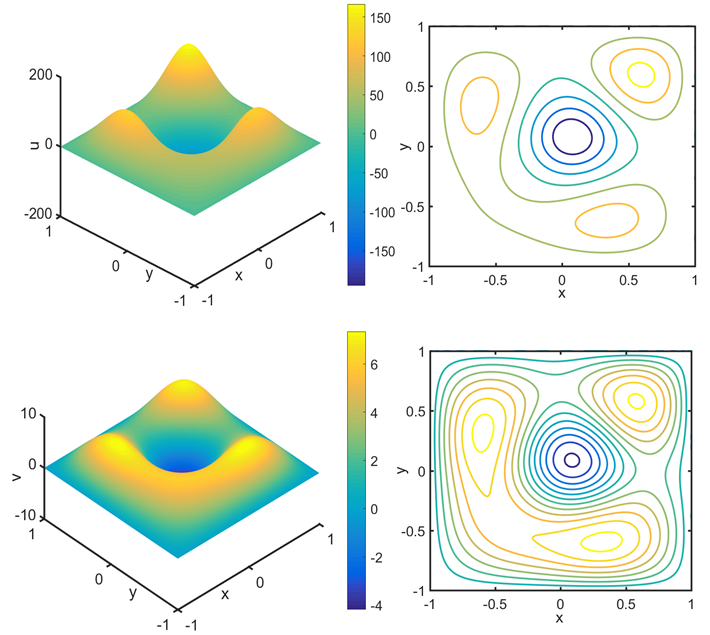}}\;
\subfigure[VI]{ \includegraphics[width=5.5cm]{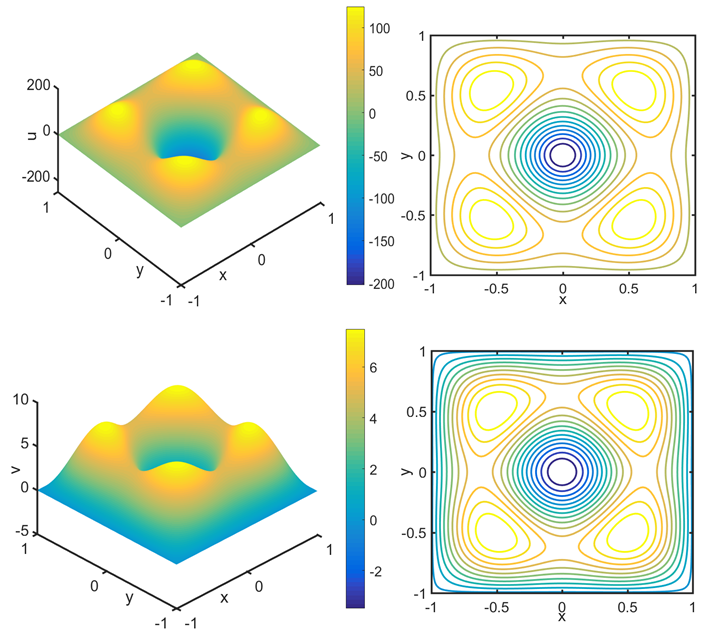}}\\
\caption{Profiles of multi-peak solutions for $p = 2, q = 2$.}
\label{pq2example4}
\end{center}
\end{figure}

When $p = 2$ and $q = 2$, from \eqref{example3} and Fig.\ref{pq2example1}, we can conclude that $G_{1}$, $G_{2} \in C^{1}[-1, 1]$ is only satisfied. As said in the introduction, the method with using second derivative information can't be used to solve \eqref{example3}. While the trust region Levenberg-Marquardt method introduced in section \ref{sect2.3} has great advantages, and numerical results are presented in Figs. \ref{pq2example2}-\ref{pq2example4} and Table. \ref{PDE1modifytable3}, indicating that our algorithm is very effective again.

\begin{table}[!h]
\centering
\caption{\small \normalsize Accuracy of our spectral trust region LM-Deflation method for $p = 2$ and $q = 2$.}
\label{PDE1modifytable3}\small
\setlength{\tabcolsep}{1.1mm}{
\begin{tabular}{c c c c c c c c}
	\hline
~   &$N$ &I  &II   &III    &IV   &V  &VI \\\cline{2-8}
\multirow{3}{*}{$||u_{24}-\hat{u}||_{\infty}$}&8 &5.45e-2 &1.79e-2 &5.93e-3 &3.74e-1 &5.86e-1 &9.32e-1\\
~  &16 &4.82e-6 &3.81e-5 &9.65e-7 &2.73e-6 &4.82e-6 &9.15e-5\\
~  &24 &3.92e-9 &8.19e-8 &2.87e-9 &9.04e-10 &8.73e-10 &6.16e-8\\\hline\hline
\multirow{3}{*}{$||v_{24}-\hat{v}||_{\infty}$}&8 &6.53e-2 &7.30e-2 &7.81e-1 &8.47e-2 &1.83e-1 &5.26e-3\\
~  &16 &4.82e-8 &7.73e-8 &6.73e-5 &6.10e-7 &1.09e-9 &5.83e-7\\
~  &24 &8.27e-10 &6.14e-12 &6.97e-11 &9.94e-9 &5.83e-11 &4.03e-10\\\hline
\end{tabular}}
\end{table}

\begin{figure}[!ht]
\begin{center}
\subfigure[I]{ \includegraphics[width=5.5cm]{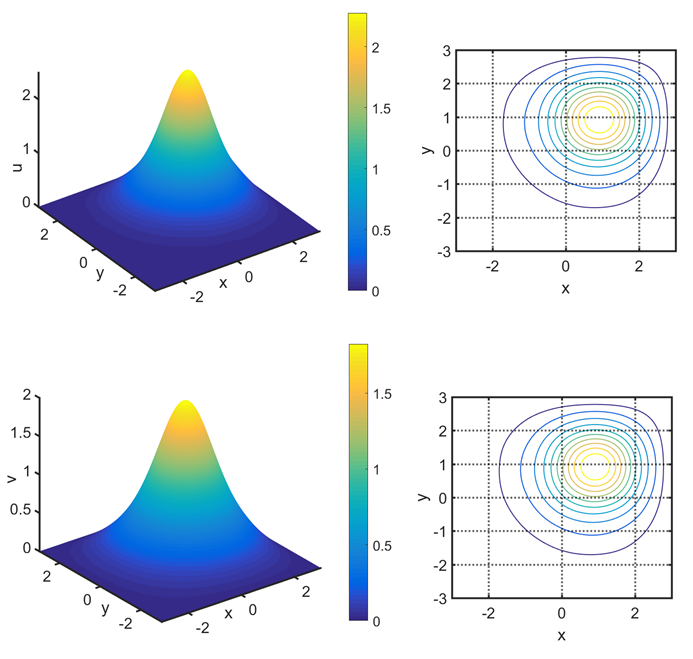}}\quad
\subfigure[II]{ \includegraphics[width=5.5cm]{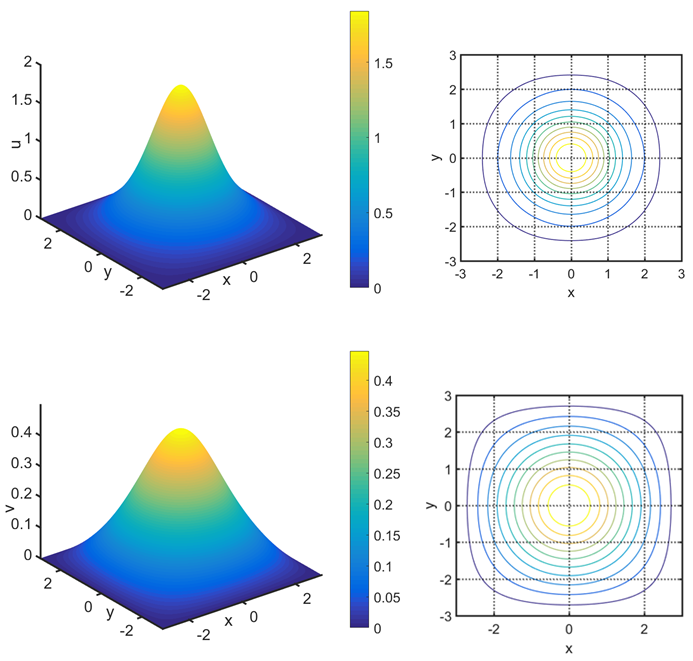}}\\\vspace{-0.4cm}
\caption{Profiles of an asymmetric positive solution (a) and a symmetric positive solution (b).}
\label{example3.21}
\end{center}
\end{figure}

\begin{figure}[!ht]
\begin{center}
\subfigure[III]{ \includegraphics[width=5.5cm]{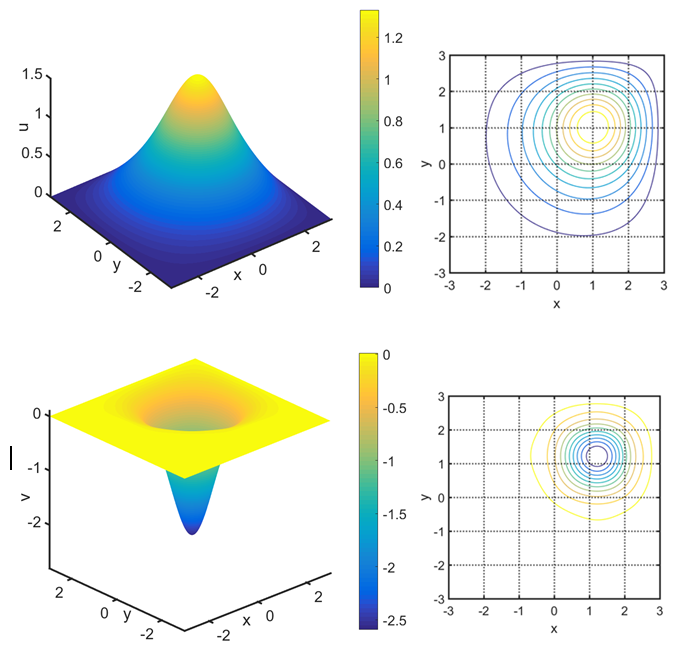}}\quad
\subfigure[IV]{ \includegraphics[width=5.5cm]{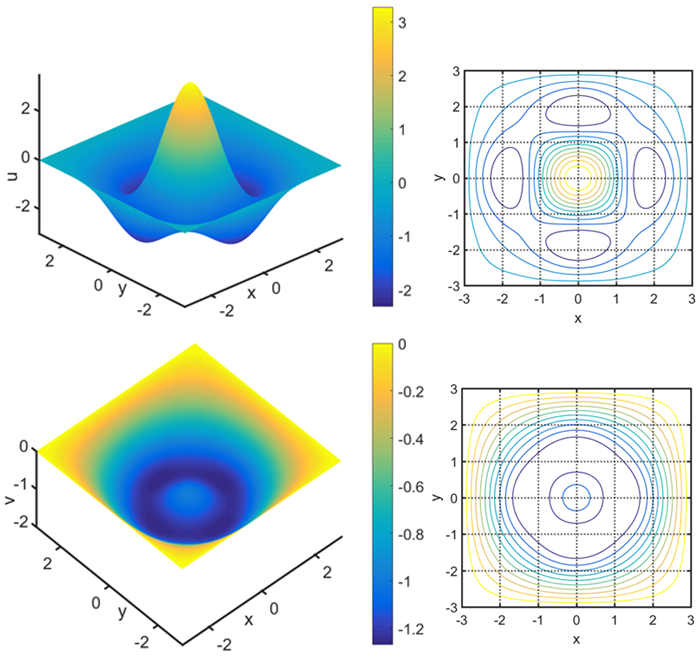}}\\\vspace{-0.4cm}
\caption{Profiles of an asymmetric positive solution (a) and a symmetric sign-changing solution (b).}
\label{example3.22}
\end{center}\vspace{-0.3cm}
\end{figure}

\begin{figure}[!ht]
\begin{center}
\subfigure[V]{ \includegraphics[width=5.5cm]{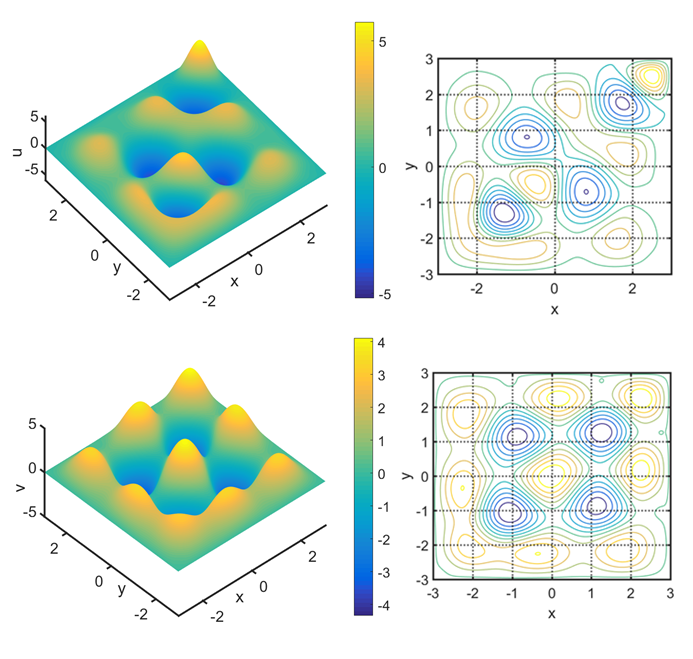}}\quad
\subfigure[VI]{ \includegraphics[width=5.5cm]{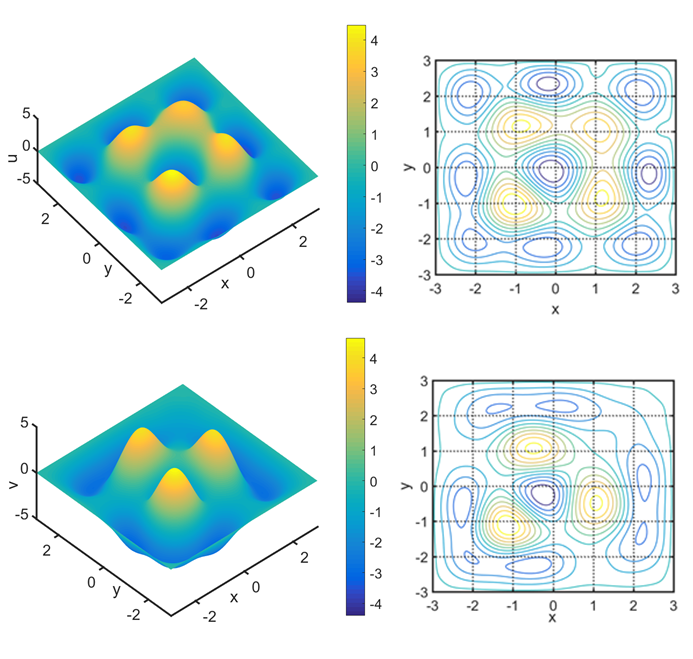}}\\\vspace{-0.4cm}
\caption{Profiles of multi-peak solutions to case 2 for $p = 3$ and $q = 3$.}
\label{example3.23}
\end{center}
\end{figure}

\noindent\underline{{\bf Case 2:}} The problem \eqref{example3} with $G_1(u(x), v(x)) = \lambda u - \delta v + |u|^{p-1}u$ and $G_2(u(x), v(x)) = \delta u + \gamma v + |v|^{q-1}v$ is known as the noncooperative system of indefinite type \cite{2010A}, which together with the domain $\Omega = (-3, 3)\times(-3, 3)$. The parameters $p = q = 3, \lambda = -0.5, \gamma = -1, \delta = 0.5$ are chosen. In Figs. \ref{example3.21}-\ref{example3.23}, several multiple solutions and their contours are presented, where some asymmetric or multi-peak solutions are also found and shown. As for the efficiency of our algorithm, the situation is similar to the case above, and we don't repeat and show it.

\begin{figure}[!ht]
\begin{center}
\subfigure[I]{ \includegraphics[width=5.5cm]{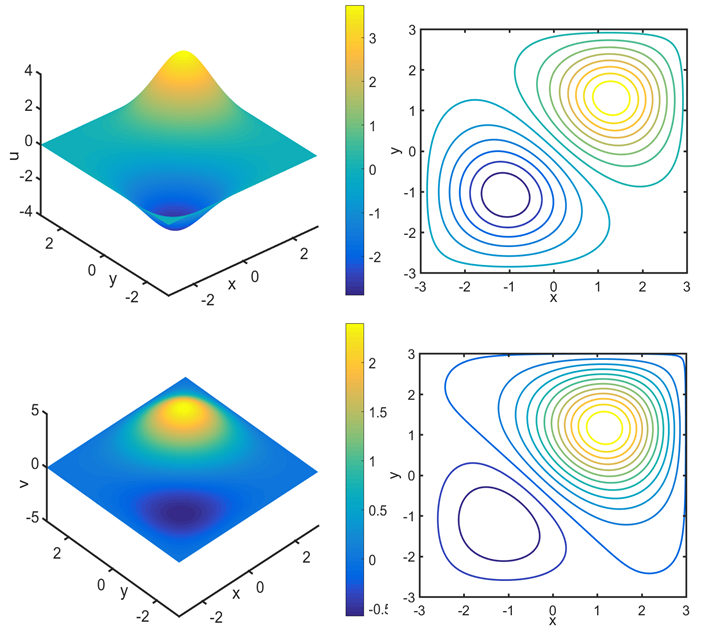}}\quad
\subfigure[II]{ \includegraphics[width=5.5cm]{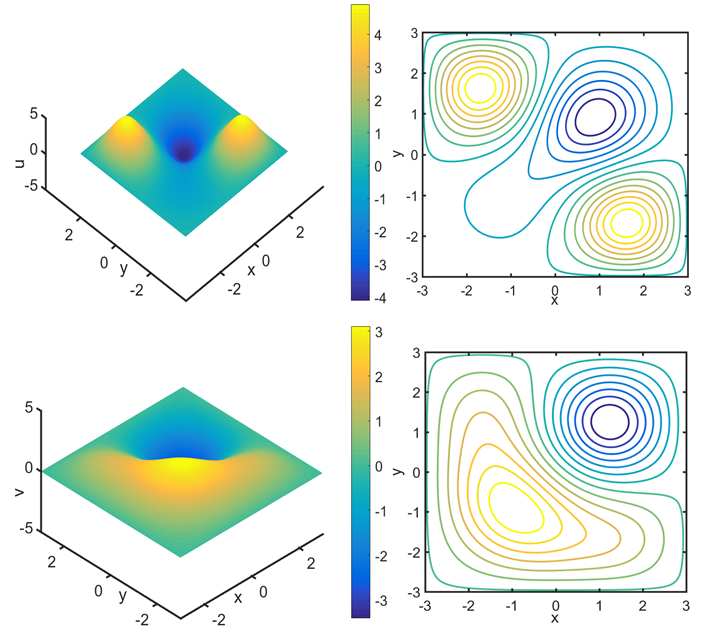}}\\
\caption{Profiles of solutions I and II for $p = 2$ and $q = 2$.}
\label{pq2example22}
\end{center}
\end{figure}

\begin{figure}[!ht]
\begin{center}
\subfigure[III]{ \includegraphics[width=5.5cm]{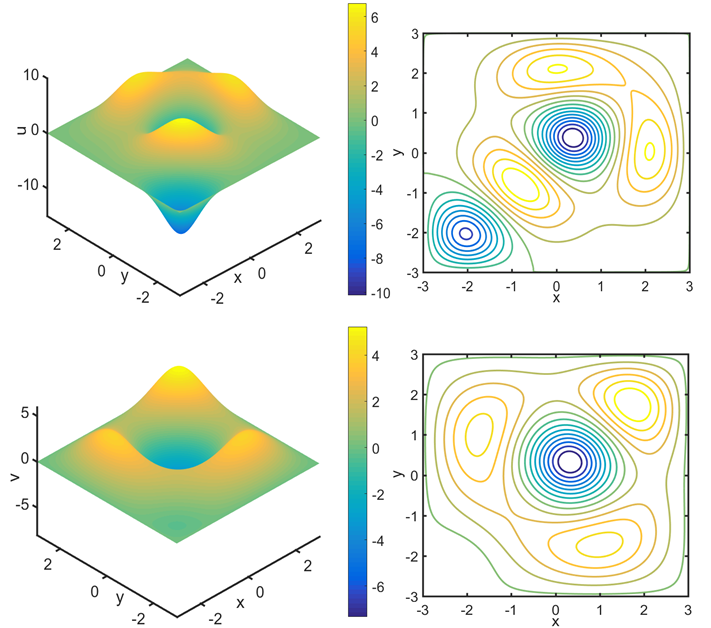}}\quad
\subfigure[IV]{ \includegraphics[width=5.5cm]{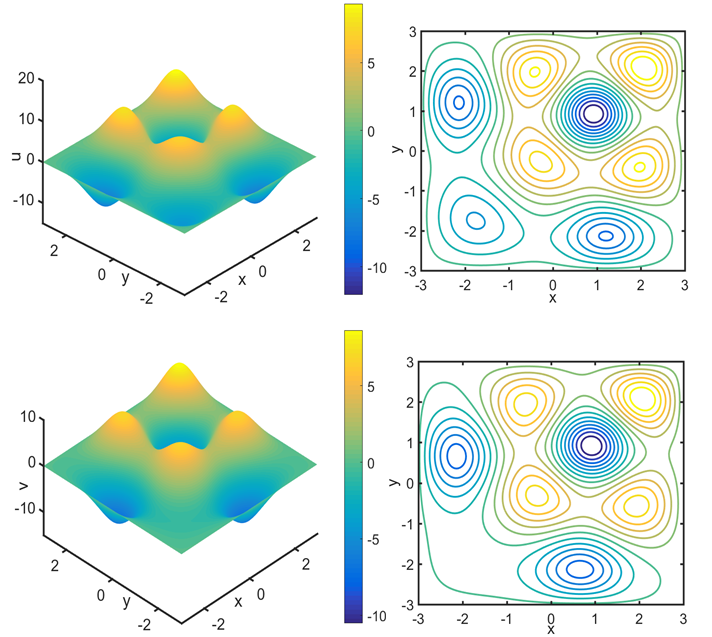}}\\
\caption{Profiles of solutions III and IV for $p = 2$ and $q = 2$.}
\label{pq2example34}
\end{center}
\end{figure}

Similar to the case 1, $G_{1}, G_{2} \in C^{1}[-1, 1]$ is only satisfied when $p = 2$ and $q = 2$. Here our aim is to test the efficiency of our algorithm, not to find multiple solutions as many as possible. Numerical results are presented in Figs.\ref{pq2example22}-\ref{pq2example34}, and the performances of the efficiency in finding multiple solutions is similar to the case 1 above. Here we don't repeat and show it.

\begin{figure}[!ht]
\begin{center}
\subfigure[I]{ \includegraphics[width=5.5cm]{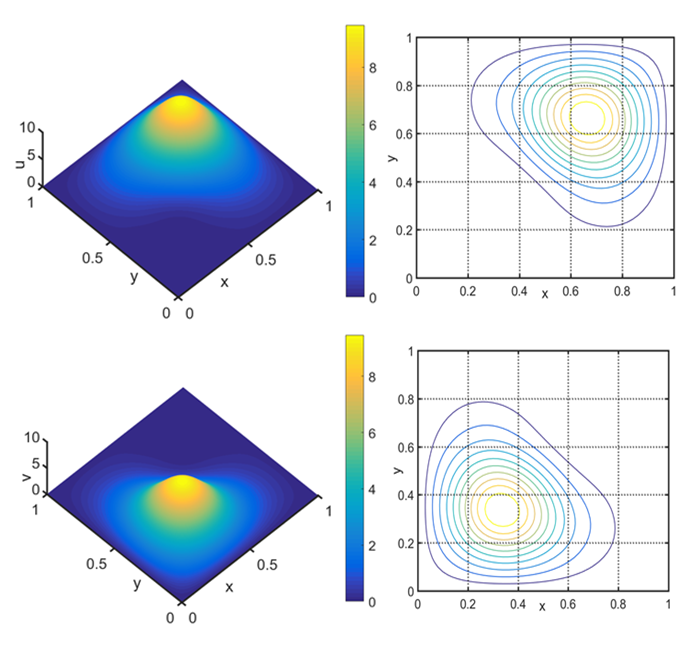}}\quad
\subfigure[II]{ \includegraphics[width=5.5cm]{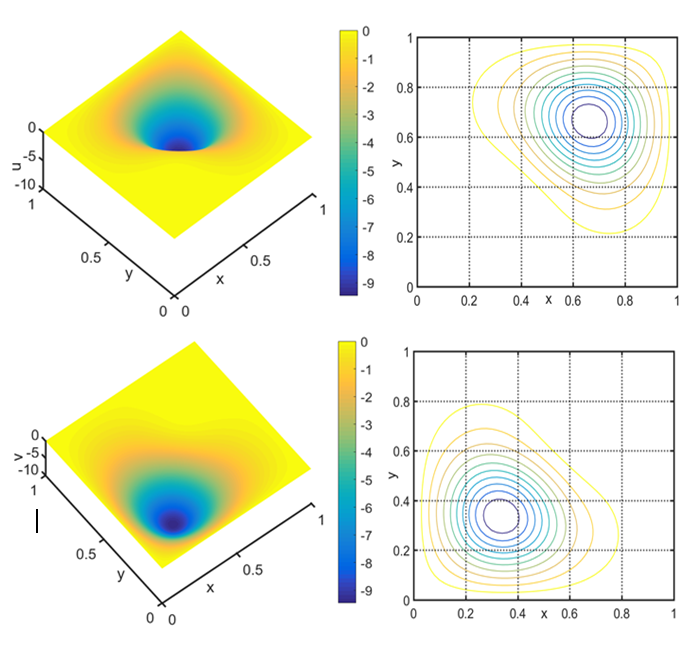}}\\
\caption{Profiles of solution set 1.}
\label{example3.31}
\end{center}
\end{figure}

\newpage

\noindent\underline{{\bf Case 3:}} The problem \eqref{example3} with $G_1(u(x), v(x)) = \eta_1 u + \mu_1 u^3 + \beta u v^2$ and $G_2(u(x), v(x)) = \eta_2 v + \mu_2 v^3 + \beta u^2v$ arises from the Bose-Einstein condensates (BEC)\cite{2020Finding}, which together with the domain $\Omega = (0, 1)\times(0, 1)$. The solution set of the problem has additional structural symmetries, that is, if

\begin{figure}[!ht]
\begin{center}
\subfigure[III]{ \includegraphics[width=5.5cm]{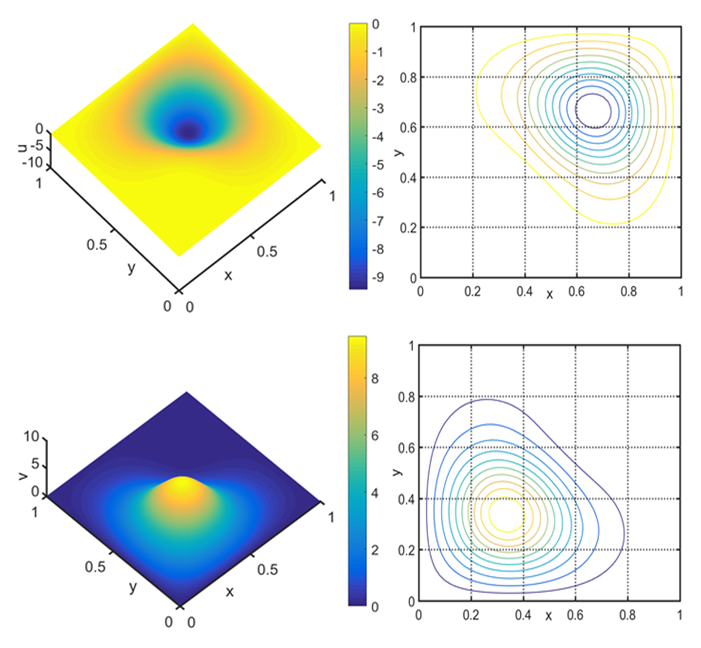}}\quad
\subfigure[IV]{ \includegraphics[width=5.5cm]{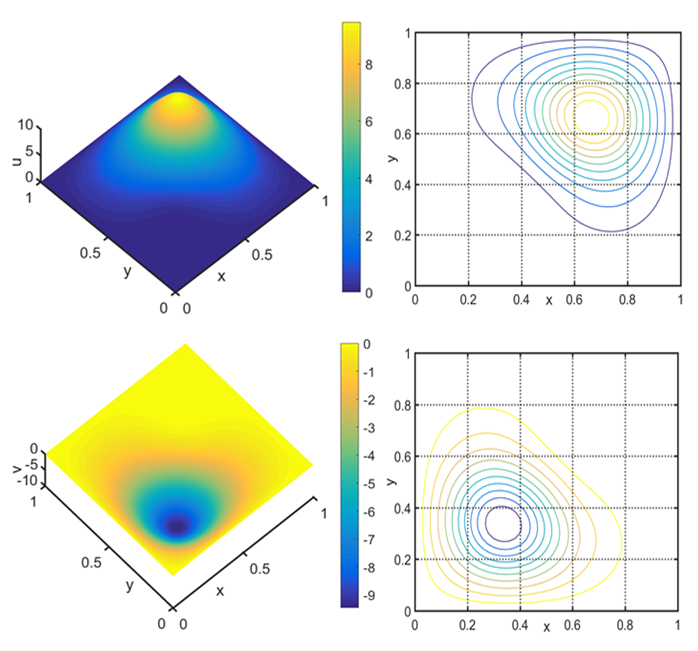}}\\
\caption{Profiles of solution set 2.}
\label{example3.32}
\end{center}\vspace{-0.6cm}
\end{figure}

\begin{figure}[!ht]
\begin{center}
\subfigure[V]{ \includegraphics[width=5.5cm]{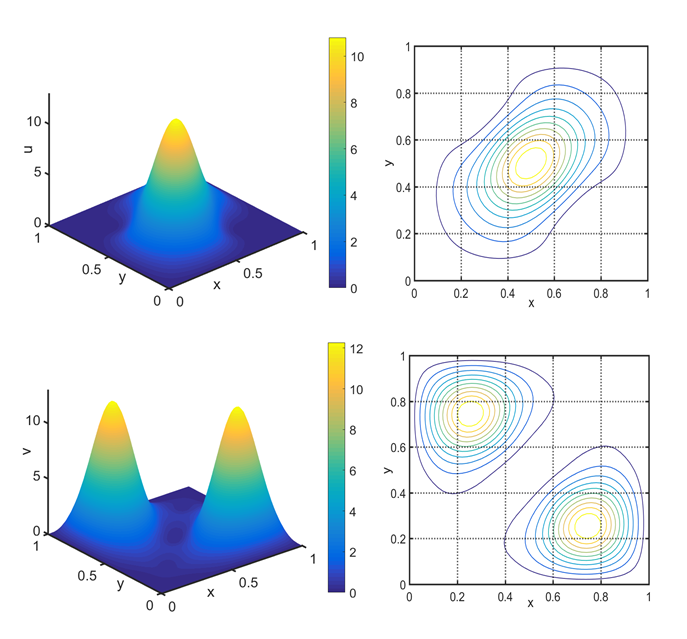}}\quad
\subfigure[VI]{ \includegraphics[width=5.5cm]{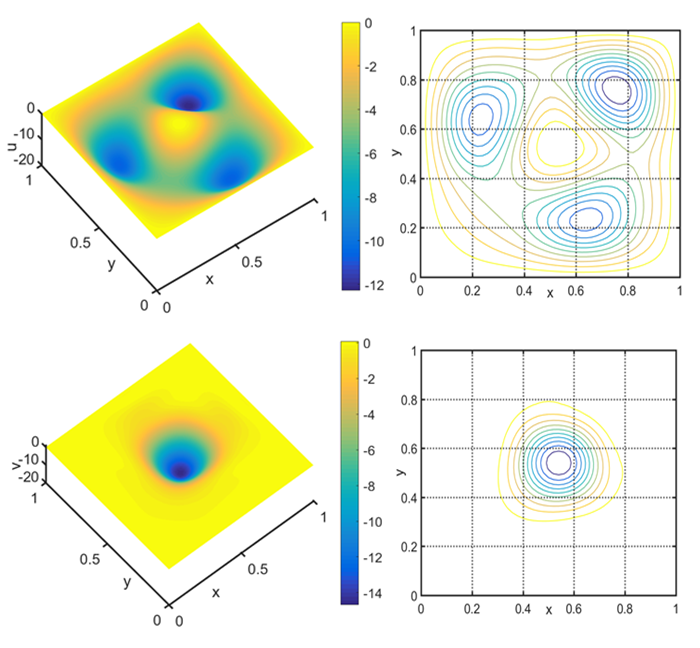}}\\
\caption{Profiles of solution set 3.}
\label{example3.33}
\end{center}\vspace{-0.3cm}
\end{figure}

\begin{figure}[!ht]
\begin{center}
\subfigure[VII]{\includegraphics[width=5.5cm]{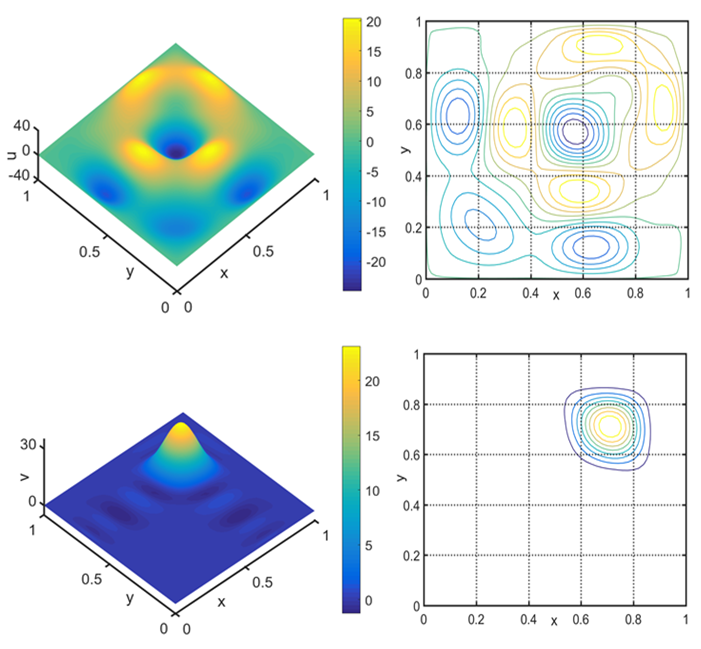}}\quad
\subfigure[VIII]{\includegraphics[width=5.5cm]{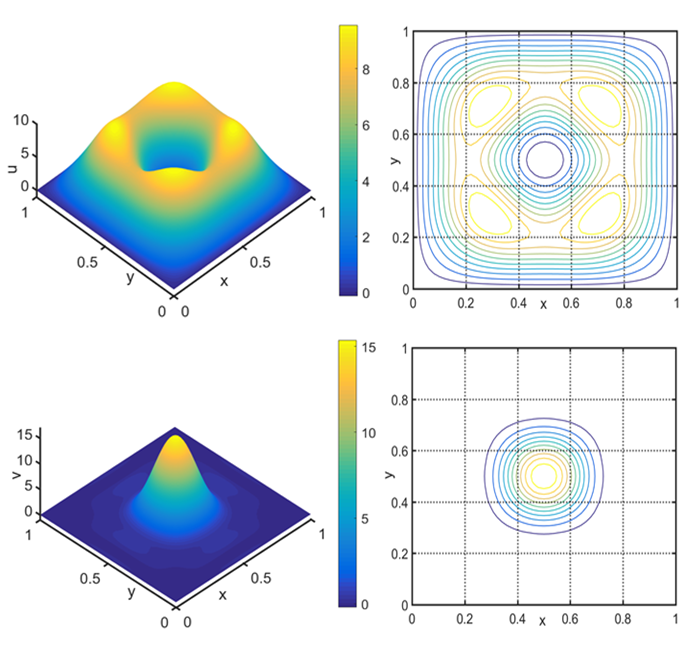}}\\
\caption{Profiles of solution set 4.}
\label{example3.34}
\end{center}
\end{figure}

\noindent $(u, v)$ is a solution, then $(\pm u, \pm v)$ are solutions. With numerical tests, we choose $\eta_1 = \eta_2 = -1, \mu_1 = \mu_2 = 1, \beta = -5.$ In Figs.\ref{example3.31}-\ref{example3.34}, we present multiple solutions, and some multiple solutions also show the symmetry. For example, if $(u, v)$ represents the type-I solution in Fig.\ref{example3.31}a, the type-II solution in Fig.\ref{example3.31}b is $(-u, -v)$. A similar situation is also observed in Fig.\ref{example3.32}. Multiple solutions presented in Figs.\ref{example3.33}-\ref{example3.34} also have the symmetry, and here we will not repeat them. In Tables \ref{PDE3modifytable1}-\ref{PDE3modifytable2} with different initial guesses in \eqref{PDE3}, numerical results on these multiple solutions are given, where the efficiency of our method is shown once again.

\begin{equation}\label{PDE3}
      IG_{1}:
      \begin{cases}
      \tilde{u}^{(0)} = \textrm{ones}(N+1, N+1)\\
      \tilde{v}^{(0)} = -\textrm{ones}(N+1, N+1)
      \end{cases}
      \;
      IG_{2}:
      \begin{cases}
      \tilde{u}^{(0)} = -\sin(\textrm{ones}(N+1, N+1))\\
      \tilde{v}^{(0)} = -\cos(\textrm{ones}(N+1, N+1)),
      \end{cases}
      \end{equation}

\begin{table}[!h]
\centering
\caption{\small \normalsize Accuracy of our spectral trust region LM-Deflation method to case 3.}
\label{PDE3modifytable1}\small
\setlength{\tabcolsep}{1.1mm}{
\begin{tabular}{c c c c c c c c c c}
	\hline
~   &$N$ &I  &II   &III    &IV   &V  &VI  &VII    &VIII \\\cline{2-10}
\multirow{3}{*}{$||u_{24}-\hat{u}||_{\infty}$}&8 &3.14e-1 &7.06e-2 &8.23e-1 &4.38e-1 &4.89e-3 &7.43e-1 &3.18e-2 &6.94e-2\\
~  &16 &3.17e-7 &3.81e-6 &4.45e-5 &3.92e-6 &2.76e-6 &7.65e-7  &4.61e-6 &7.95e-6\\
~  &24 &6.55e-9 &5.52e-9 &1.86e-9 &7.09e-11 &7.54e-10 &9.73e-9  &2.76e-11 &6.94e-10 \\\hline\hline
\multirow{3}{*}{$||v_{24}-\hat{v}||_{\infty}$}&8 &7.13e-1 &5.05e-2 &2.54e-1 &8.14e-2 &9.59e-2 &6.79e-2  &2.55e-1 &8.40e-1\\
~  &16 &2.76e-7 &1.62e-7 &1.19e-6 &6.99e-7 &1.49e-6 &5.47e-6  &6.79e-6 &5.85e-5 \\
~  &24 &3.80e-11 &7.79e-11 &5.30e-12 &2.51e-11 &7.53e-11 &1.19e-11  &7.58e-10 &4.69e-10\\\hline
\end{tabular}}
\end{table}

\begin{table}[!h]
\centering\small
\caption{\small A comparison of $\|\textbf{\textit{F}}(\textbf{\textit{x}})\|_{2}$ between our method and other methods for case 3.}
\label{PDE3modifytable2}
\setlength{\tabcolsep}{1.1mm}{
\begin{tabular}{ccc|ccc|ccc}
	\hline
\multicolumn{3}{c}{Newtonian iteration in \cite{farrell2015deflation}} &\multicolumn{3}{c}{LSTR method in \cite{2022Two}} &\multicolumn{3}{c}{our method}\\\hline
$n_{it}$ & $IG_{1}$  & $IG_{2}$   &$n_{it}$ & $IG_{1}$  & $IG_{2}$  &$n_{it}$ & $IG_{1}$  & $IG_{2}$\\
5 &1.02e8 &6.01e10          &5    &4.50e6  &8.25e8     &5    &4.42e8  &1.52e9 \\
10  &7.94e15  &2.62e17      &10   &1.06e3   &9.61e3    &10   &9.13e4  &7.79e5 \\
15  &3.11e23  &6.98e25      &15   &6.54e1   &2.28e0    &15   &6.89e1  &6.54e0 \\
20  &7.01e27  &9.13e35      &20   &9.96e-3  &7.81e-3   &20   &7.81e-5  &5.28e-6 \\
~  &- &-                    &25   &5.28e-7 &2.28e-7    &25   &4.63e-10 &1.65e-11 \\\hline
T &- &-                     &~    &8.01     &9.21      &~     &6.09    &5.97   \\\hline
\end{tabular}}
\end{table}

\section{Concluding remarks}\label{sect4}

In this paper, an efficient trust region LM-Deflation method was proposed to find multiple solutions of \eqref{New1.1}. We demonstrated that the trust region LM-Deflation method was very effective. In subsequent papers, we will continue to study this method including its convergence analysis.

\bibliographystyle{siam}
\bibliography{myreference1}

\begin{thebibliography}{10}

\bibitem{allgower2009application}
E.~L. Allgower, S.~G. Cruceanu, and S.~Tavener.
\newblock Application of numerical continuation to compute all solutions of
  semilinear elliptic equations.
\newblock {\em Adv. Geom}, 76(2009):pp. 1--10.

\bibitem{allgower2006solution}
E.~L. Allgower, A.~J.~Sommese D.~J.~Bates, and C.~W. Wampler.
\newblock Solution of polynomial systems derived from differential equations.
\newblock {\em Computing}, 76(2006):pp. 1--10.

\bibitem{1971Deflation}
K.~M. Brow and W.~B. Gearhart.
\newblock Deflation techniques for the calculation of further solutions of a
  nonlinear system.
\newblock {\em Numer. Math.}, 16(1971):334--342.

\bibitem{chen2004search}
C.~M. Chen and Z.~Q. Xie.
\newblock Search extension method for multiple solutions of a nonlinear
  problem.
\newblock {\em Comp. Math. Appl.}, 47(2004):pp. 327--343.

\bibitem{2010A}
X.~J. Chen and J.~X. Zhou.
\newblock A local min-max-orthogonal method for finding multiple solutions to
  noncooperative elliptic systems.
\newblock {\em Math. Comp.}, 79(2010):pp. 2213--2236.

\bibitem{choi1993mountain}
Y.~S. Choi and P.~J. McKenna.
\newblock A mountain pass method for the numerical solution of semilinear
  elliptic problems.
\newblock {\em Nonlinear. Anal.}, 20(1993):pp. 417--437.

\bibitem{davis1960introduction}
H.~T. Davis.
\newblock {\em Introduction to nonlinear differential and integral equations}.
\newblock US Atomic Energy Commission, 1960.

\bibitem{ding1999high}
Z.~H. Ding, D.~Costa, and G.~Chen.
\newblock A high-linking algorithm for sign-changing solutions of semilinear
  elliptic equations.
\newblock {\em Nonlinear. Anal.}, 38 (1999):pp. 151--172.

\bibitem{2012Fan}
J.~Y. Fan.
\newblock The modified levenberg-marquardt method for nonlinear equations with
  cubic convergence.
\newblock {\em Math. Comp}, 81(2012):pp. 447--466.

\bibitem{2019Fan}
J.~Y. Fan, J.~C. Huang, and J.~Y. Pan.
\newblock An adaptive multi-step levenberg–marquardt method.
\newblock {\em J. Sci. Comput}, 78(2019):pp. 531--548.

\bibitem{2005Fan}
J.~Y. Fan and Y.~X. Yuan.
\newblock On the quadratic convergence of the levenberg-marquardt method
  without nonsingularity assumption.
\newblock {\em Computing}, 74(2019):pp. 23--39.

\bibitem{farrell2015deflation}
P.~E. Farrell, A.~Birkisson, and S.~W. Funke.
\newblock Deflation techniques for finding distinct solutions of nonlinear
  partial differential equations.
\newblock {\em SIAM J. Sci. Comput.}, 37(2015):pp. A2026--A2045.

\bibitem{2002A}
U.~Frisch, S.~Matarrese, R.~Mohayaee, and A.~Sobolevski.
\newblock A reconstruction of the initial conditions of the universe by optimal
  mass transportation.
\newblock {\em Nature}, 417(2002):pp. 260--262.

\bibitem{hao2014bootstrapping}
W.~R. Hao, J.~D. Hauenstein, B.~Hu, and A.~J. Sommese.
\newblock A bootstrapping approach for computing multiple solutions of
  differential equations.
\newblock {\em J. Comput. Appl. Math.}, 258(2014):pp. 181--190.

\bibitem{2020Spatial}
W.~R. Hao and C.~Xue.
\newblock Spatial pattern formation in reaction–diffusion models: a
  computational approach.
\newblock {\em J. Math. Biol.}, 80(2020):pp. 521--543.

\bibitem{2022Two}
L.~Li, L.~L. Wang, and H.~Y. Li.
\newblock An efficient spectral trust-region deflation method for multiple
  solutions.
\newblock {\em J. Sci. Comput}, 32(2023):pp. 1--23.

\bibitem{li2001minimax}
Y.~X. Li and J.~X. Zhou.
\newblock A minimax method for finding multiple critical points and its
  applications to semilinear pdes.
\newblock {\em SIAM J. Sci. Comput.}, 23(2001):pp. 840--865.

\bibitem{shen2011spectral}
J.~Shen, T.~Tang, and L.~L. Wang.
\newblock {\em Spectral methods: algorithms, analysis and applications},
  volume~41.
\newblock Springer Science \& Business Media, 2011.

\bibitem{sun2006optimization}
W.~Y. Sun and Y.~X. Yuan.
\newblock {\em Optimization theory and methods: nonlinear programming},
  volume~1.
\newblock Springer Science \& Business Media, 2006.

\bibitem{Tadmor2012A}
E.~Tadmor.
\newblock A review of numerical methods for nonlinear partial differential
  equations.
\newblock {\em B. Am. Math. Soc.}, 49(2012):pp. 507--554.

\bibitem{2018Two}
Y.~Wang, W.~Hao, and G.~Lin.
\newblock Two-level spectral methods for nonlinear elliptic equations with
  multiple solutions.
\newblock {\em SIAM J. Sci. Comput}, 40(2018):pp. B1180--B1205.

\bibitem{XIE2025116146}
P.~Xie.
\newblock Sufficient conditions for error distance reduction in the
  \(\ell^2\)-norm trust region between minimizers of local nonconvex
  multivariate quadratic approximates.
\newblock {\em Journal of Computational and Applied Mathematics}, 453:116146,
  2025.

\bibitem{xie2025remuregionalminimalupdating}
P.~Xie and S.~M. Wild.
\newblock {ReMU}: Regional minimal updating for model-based derivative-free
  optimization.
\newblock 2025.

\bibitem{xie2023derivative}
P.~Xie and Y.~Yuan.
\newblock A derivative-free optimization algorithm combining line-search and
  trust-region techniques.
\newblock {\em Chinese Annals of Mathematics, Series B}, 44(5):719--734, 2023.

\bibitem{xie2024derivative}
P.~Xie and Y.~Yuan.
\newblock Derivative-free optimization with transformed objective functions and
  the algorithm based on the least {Frobenius} norm updating quadratic model.
\newblock {\em Journal of the Operations Research Society of China}, pages
  1--37, 2024.

\bibitem{xie2024newtwodimensionalmodelbasedsubspace}
P.~Xie and Y.~Yuan.
\newblock A new two-dimensional model-based subspace method for large-scale
  unconstrained derivative-free optimization: {2D-MoSub}, 2024.

\bibitem{xie2023least}
P.~Xie and Y.~X. Yuan.
\newblock Least {$H^2$} norm updating of quadratic interpolation models for
  derivative-free trust-region algorithms.
\newblock {\em IMA Journal of Numerical Analysis}, page drae106, 03 2025.

\bibitem{xie2005improved}
Z.~Q. Xie, C.~M. Chen, and Y.~Xu.
\newblock An improved search-extension method for computing multiple solutions
  of semilinear {PDE}s.
\newblock {\em IMA J. Numer. Anal.}, 25(2005):pp. 549--576.

\bibitem{xie2006improved}
Z.~Q. Xie, C.~M. Chen, and Y.~Xu.
\newblock An improved search-extention method for solving semilinear {PDE}s.
\newblock {\em Acta. Math. Sci.}, 26(2006):pp. 757--766.

\bibitem{xie2015augmented}
Z.~Q. Xie, W.~F. Yi, and J.~X. Zhou.
\newblock An augmented singular transform and its partial newton method for
  finding new solutions.
\newblock {\em J. Comput. Appl. Math.}, 286(2015):pp. 145--157.

\bibitem{2001Yamashita}
N.~Yamashita and M.~Fukushima.
\newblock {\em On the rate of convergence of the Levenberg-Marquardt method},
  volume~15.
\newblock Computing (Supplement 15), 2001.

\bibitem{2005A}
X.~D. Yao and J.~X. Zhou.
\newblock A minimax method for finding multiple critical points in banach
  spaces and its application to quasi-linear elliptic {PDE}s.
\newblock {\em SIAM J. Sci. Comput.}, 26(2005):pp. 1796--1809.

\bibitem{2007numerical}
X.~D. Yao and J.~X. Zhou.
\newblock Numerical methods for computing nonlinear eigenpairs: Part {I}.
  {I}so-{H}omogeneous cases.
\newblock {\em SIAM J. Sci. Comput.}, 29(2007):pp. 1355--1374.

\bibitem{2008numerical}
X.~D. Yao and J.~X. Zhou.
\newblock Numerical methods for computing nonlinear eigenpairs: Part {II}.
  {N}on-{I}so-{H}omogeneous cases.
\newblock {\em SIAM J. Sci. Comput}, 30(2008):pp. 937--956.

\bibitem{2020Finding}
X.~P. Zhang, J.~Zhang, and B.~Yu.
\newblock Finding multiple solutions to elliptic systems with polynomial
  nonlinearity.
\newblock {\em Numer. Meth. Part. D. E.}, 36(2019):pp. 1074--1097.

\bibitem{zhang2013eigenfunction}
X.~P. Zhang, J.~T Zhang, and B.~Yu.
\newblock Eigenfunction expansion method for multiple solutions of semilinear
  elliptic equations with polynomial nonlinearity.
\newblock {\em SIAM J. Numer. Anal.}, 51(2013):pp. 2680--2699.

\bibitem{2012ZZhang}
Z.~K. Zhang.
\newblock The research on derivative-free optimization methods.
\newblock {\em PhD thesis, University of Chinese Academy of Sciences}, 2012.

\bibitem{2014Zhangzaikun}
Z.~K. Zhang.
\newblock Sobolev seminorm of quadratic functions with applications to
  derivative-free optimization.
\newblock {\em Math. Program}, 146(2014):pp. 77-- 96.

\bibitem{2005Instability}
J.~X. Zhou.
\newblock Instability analysis of saddle points by a local minimax method.
\newblock {\em Math. Comp.}, 74(2004):pp. 1391--1411.

\end{thebibliography}

\end{document}